\documentclass[10pt,a4paper]{article}
\pdfoutput=1
\usepackage{amsmath,amsthm, amssymb, url, enumerate}
\usepackage{graphicx}
\usepackage{cases}
\usepackage{epic, eepic, ecltree}
\usepackage{stmaryrd}
\usepackage[T1]{fontenc}
\usepackage{textcomp}
\usepackage{color}
\usepackage[all]{xy}
\usepackage{cases}
\usepackage{stmaryrd}
\usepackage{bm}
\usepackage{here}

\theoremstyle{plain}
\newtheorem{defi}{Definition}[section]
\newtheorem{thm}[defi]{Theorem}
\newtheorem{prop}[defi]{Proposition} 
\newtheorem{lem}[defi]{Lemma}
\newtheorem{rem}[defi]{Remark}
\newtheorem{cor}[defi]{Corollary}
\newtheorem{conj}[defi]{Conjecture}

\newtheorem{ex}[defi]{Example}

\renewcommand{\proofname}{Proof.\ }
\renewcommand{\thefootnote}{\fnsymbol{footnote}}
\renewcommand{\thefootnote}{} 
\newcommand{\slmc}[1]{SL(#1,\mathbb{C})}
\newcommand{\glmc}[1]{GL(#1,\mathbb{C})}

\newcommand{\compl}[1]{\mathbb{C}^{#1}}
\newcommand{\real}[1]{\mathbb{R}^{#1}}
\newcommand{\ratio}[1]{\mathbb{Q}^{#1}}
\newcommand{\inte}[1]{\mathbb{Z}^{#1}}
\newcommand{\natu}[1]{\mathbb{N}^{#1}}

\newcommand{\age}{\mathrm{age}}
\newcommand{\AHilb}{\mathrm{A\mathchar`-Hilb}(\mathbb{C}^3)}
\newcommand{\CFF}{\mathrm{CFF}}
\newcommand{\GHilb}[1]{G\mathrm{\mathchar`-Hilb}(\mathbb{C}^{#1})}
\newcommand{\Hilb}{\mathrm{Hilb}}
\newcommand{\rank}{\mathrm{rank}}

\newcommand{\lcm}{\mathrm{LCM}}
\newcommand{\Ht}{\mathrm{ht}}
\newcommand{\Conv}{\mathrm{Conv}}

\makeatletter
\def\mapstofill@{%
   \arrowfill@{\mapstochar\relbar}\relbar\rightarrow}
\newcommand*\xmapsto[2][]{%
   \ext@arrow 0395\mapstofill@{#1}{#2}}
\makeatother

\begin{document}
\title{\vspace{-3cm}Crepant Property of Fujiki-Oka Resolutions for\\ Gorenstein Abelian Quotient Singularities}
\author{Kohei\ Sato \quad Yusuke\ Sato}
\date{}
\maketitle 
\thispagestyle{empty}
\footnote{ 
2010 \textit{Mathematics Subject Classification}.
Primary 14B05, 14J17; Secondary 13H10, 14C17, 14J30, 14J35, 14J40, 14L30, 14M25, 52B20.
}
\footnote{ 
\textit{Key words and phrases}.
Crepant resolutions, Fujiki-Oka resolutions, Multidimensional continued fractions, Hirzebruch-Jung continued fractions, Quotient singularities, Higher dimension, Finite groups, Abelian groups, Toric varieties, Invariant theory.
}
\begin{abstract}
We show a necessary and sufficient condition for the Fujiki-Oka resolutions of Gorenstein abelian quotient singularities to be crepant in all dimensions by using Ashikaga's continuous fractions. Moreover, we prove that all three dimensional Gorenstein abelian quotient singularities always possess a crepant Fujiki-Oka resolution as a corollary. This alternative proof of existence needs only simple computations comparing with the results ever known.
\end{abstract}

\section{Introduction}
\noindent

Let $G$ be a finite subgroup of $\glmc{n}$, and let $(\compl{n}/G, [0])$ be an $n$-dimensional quotient singularity. In the case $n=2$, the cyclic quotient singularities $(\compl{2}/G, [0])$ have an unique minimal resolution, and the self-intersection numbers of the exceptional divisors of the minimal resolution correspond to the coefficients of the Hirzebruch-Jung continued fraction related to the group action of $G$ (See Section \ref{HJCF}). Moreover, if $G$ is a finite subgroup of $\slmc{2}$, then the quotient singularity is Gorenstein, and the dual of the weighted graph obtained from the exceptional divisors of the minimal resolution corresponds to the Dynkin diagram obtained from the non-trivial irreducible representations of $G$. This correspondence remarked by J. McKay \cite{25} is called the{\it McKay correspondence}. 

This correspondence has not been shown in general dimension. It is because minimal resolutions do not necessarily exist for cyclic quotient singularities in the case of $n\geq 3$. On the other hand, the McKay correspondence has been generalized to the case of $n=3$ by V. B. Batyrev and D. I. Dais \cite{2} and by Y.Ito and M. Reid \cite{16} as the following; 
\[ \{{\rm Conjugacy\ classes\ of}\ G\ {\rm of}\ \age\ i \} \longleftrightarrow \{ {\rm A\ basis\ of}\ H^{2i}(\widetilde{\compl{3}/G},\ratio{}) \}\]
where $\widetilde{\compl{3}/G}$ is a {\it crepant resolution} of $\compl{3}/G$, i.e., a resolution such that the canonical divisor of $\widetilde{\compl{3}/G}$ is trivial.

Therefore the existence of crepant resolutions is a necessary condition to construct the above correspondences, but, unfortunately, there does not necessarily exist a crepant resolution for arbitrary Gorenstein quotient singularities. Hence, ``The Existence Problem of Crepant Resolutions'' is a basic question to generalize the McKay correspondence.

As for known results with respect to the existence problem, if $n=2$, then minimal resolutions are crepant. In the case of $n=3$, all Gorenstein quotient singularities possess crepant resolutions. This result was proved by Y.Ito \cite{13, 14}, D. G. Markushevich \cite{20,21,23} and S. S. Roan \cite{30} case by case based on the classification of three dimensional Gorenstein quotient singularities given by S. S. Yau and Y. Yu \cite{32}. In the case of $n\geq4$, Gorenstein quotient singularities do not necessarily have a crepant resolution. On the other hand, D. I. Dais, M. Henk, G. M. Ziegler and others have proved that all the complete intersection Gorenstein quotient singularities possess a crepant resolution and have constructed some infinite series of Gorenstein quotient singularities which possess a crepant resolution \cite{5,6,7,8}. 
Moreover, some infinite series of Gorenstein quotient singularities which possess a crepant resolution were constructed by others \cite{10,31}. 

In this paper, the {\it $G$-Hilbert scheme} $\GHilb{n}$ for a finite subgroup $G$ of $\glmc{n}$ means the irreducible component in the Hilbert scheme $\mathrm{Hilb}^{r}(\compl{n})$ of $r$ points on $\compl{n}$ dominating $\compl{n}/G$ where $r$ is the order of $G$, i.e., $
\GHilb{n}:=\{ I\subset \compl{}[x_1,x_2,\ldots , x_n] | I:\text{a $G$-invariant ideal },\ \dim(\mathcal{O}_{\compl{n}}/I)=|G| \}$.
 Some relations between crepant resolutions and $G$-Hilbert schemes are known. In the case of $n=2$ (resp. $n=3$), $\GHilb{n}$ is the minimal resolution (resp. a crepant resolution) of $\compl{n}/G$ for arbitrary finite subgroup $G$ in $\glmc{2}$ (resp.$\slmc{3}$) \cite{3,11,15,18,26}. However, in the case of $n\geq4$, $\GHilb{n}$ is not necessarily smooth, and flops of crepant resolutions do not keep smoothness \cite{24}.

In this paper, we shall show a necessary and sufficient condition for Gorenstein abelian quotient singularities in arbitrary dimension to admit a crepant Fujiki-Oka resolution by using {\it Ashikaga's continued fractions} \cite{1}. The Ashikaga's continued fraction consists of a {\it remainder polynomial} and a {\it round down polynomial} (See Definition\ref{DOACF}), and the round down polynomial is a dimensional generalization of the Hirzebrch-Jung continuous fraction. The remainder polynomial (resp. The round up polynomial) indicates the types of quotient singularities (resp. the $\inte{n-1}$-weight) which appear in each step of the {\it Fujiki-Oka resolution} \cite{9, 28}. The Fjiki-Oka resolutions are always smooth. In order to prove a Fujiki-Oka resolution is crepant, it is enough to show that the Fujiki-Oka resolution preserves the crepant property at each step. In Chapter 3, we shall show this condition can be expressed by the coefficients of the remainder polynomials as follows.\\

\noindent
\textbf{Theorem \ref{thm1}.}\ 
{\it For a cyclic quotient singularity of $\frac{1}{r}(1,a_2,\dots,a_n)$-type, the Fujiki-Oka resolution is crepant if and only if the ages of all the coefficients of the corresponding remainder polynomial $\mathcal{R}_*\left(\frac{(1,a_2,\dots,a_n)}{r}\right)$ are $1$.}\\

In Chapter \ref{Ab}, we introduce an extension of the Fujiki-Oka resolutions to abelian case, which is named {\it iterated Fujiki-Oka resolutions}. By using them, we shall generalize this theorem to abelian case.\\

\noindent
\textbf{Theorem \ref{thm2}.}\ 
{\it Let $\widetilde{Y_{H_1}}, \widetilde{Y_{H_2}}, \ldots, \widetilde{Y_{H_k}}=\widetilde{Y_G}$ be a sequence of iterated Fujiki-Oka resolutions for an $n$-dimensional Gorenstein abelian quotient singularity $\compl{n}/G$. The iterated Fujiki-Oka resolution $\widetilde{Y_G}$ is crepant resolution for $\compl{n}/G$ if and only if the ages of all the coefficients in the remainder polynomials associated with every $\widetilde{Y_{H_i}}\ (i=1,\ldots, k)$ are $1$. }\\

As a corollary of this theorem, we have the following result of the existence of crepant resolutions for three dimensional Gorenstein abelian quotient singularities.\\

\noindent
\textbf{Corollary \ref{cor2}.}\ 
{\it All three dimensional Gorenstein abelian quotient singularities possess a crepant iterated Fujiki-Oka resolution.}\\

The proof of this corollary is an alternative proof of existence of crepant resolutions for abelian Gorenstein quotient singularities in dimension three, and that needs only simple computations comparing with the results ever known. In the Chapter \ref{Discuss}, we see some relations between Fujiki-Oka resolutions and $\mathrm{A}$-Hilbert schemes.

\section*{Acknowledgments.}
We would like to thank Professor Tadashi Ashikaga for dedicated support, especially, helpful discussions at Tohokugakuin University. We also thank Professor Yukari Ito for giving us many useful advices. National Institute of Technology, Oyama college has supported our study.

\section{Toric geometry and Continued fractions}
The purpose of this section is to introduce some basic notions to prove the main result. 
Let $G$ be a finite subgroup of $GL(n,\compl{})$ of order $r$. If $G$ is abelian, then all elements in $G$ are simultaneously diagonalizable. Therefore, any element in $G$ can be written as the form
$g={\rm diag}(e^{\frac{2a_1\pi \sqrt{-1}}{r}},\dots,e^{\frac{2a_n\pi \sqrt{-1}}{r}})$ where $1\leq i \leq n$ and $0\leq a_i < r$.  For simplicity, the matrix ${\rm diag}(e^{\frac{2a_1\pi \sqrt{-1}}{r}},\dots,e^{\frac{2a_n\pi \sqrt{-1}}{r}})$ is denoted by $\frac{1}{r}(a_1,\dots,a_n)$. In the following, we assume that $G$ is abelian.

\subsection{Notations from Toric Geometry}\label{TG}
Let $N$ be a free $\inte{}$-module of rank $n$ and $N_\real{} = N \otimes_{\inte{}}\real{}$. Let $\bm{e}_1,\ldots, \bm{e}_n$ be a fixed  basis of $N$. If the convex hull $\Conv\{ \bm{0},\bm{n}\}$ contains no elements in $N$ except $\bm{0}$ and $\bm{n}$, the element $\bm{n}\in N$ is called {\it primitive}. For $\bm{n}_1,\ldots ,\bm{n}_k \in N$, the subset $\tau = \real{}_{\geq0}\bm{n}_1 + \cdots + \real{}_{\geq0}\bm{n}_k \subset N_{\real{}}$ satisfying $\tau\cap (-\tau)=\bm{0}$ is called a {\it rational strongly convex polyhedral cone} where $\real{}_{\geq0}$ is the set of all non negative elements in $\real{}$. For simplicity, $\tau$ also signifies the finite fan consists of all faces of $\tau$. The {\it dimension} of a cone $\tau$ is defined as the dimension of $\real{}\cdot\tau$ as vector space over $\real{}$. If the dimension of a cone $\tau$ is $n$, then the cone is called {\it maximal}. Let $\sigma = \real{}_{\geq0}\bm{e}_1 + \cdots + \real{}_{\geq0}\bm{e}_n \subset N_{\real{}}$. The toric variety $X(N,\sigma)$ determined by $N$ and the finite fan $\sigma$ is isomorphic to $\mathbb{C}^n$. There exists a morphism of toric varieties $\phi_T : X(N,\sigma) \to X(N^{\prime},\sigma)$ corresponding to the quotient map $\phi: \mathbb{C}^n \to \mathbb{C}^n/G$ where $N^{\prime}$ is the free $\inte{}$-module of rank $n$ satisfying $N\subset N^{\prime}$ and $N^{\prime}/N \cong G$ as groups. Therefore, there is an element $\bar{g}=\frac{1}{r}(a_1,\dots,a_n) \in N^{\prime}$ for each $g\in G$. We set $N^{\prime}$ as the following:
\[ \displaystyle N^{\prime}=N + \sum_{\frac{1}{r}(a_1,\ldots,a_n) \in G}\ \frac{1}{r}(a_1,\ldots,a_n) \inte{},\]
and $\bar{g}=\frac{1}{r}(a_1,\dots,a_n) \in N^{\prime}$ maps to $g=\frac{1}{r}(a_1,\dots,a_n) \in G$ by the composition of the quotient map and the isomorphism from $N^{\prime}$ to $G$. We note that $N^{\prime}_{\real{}}$ satisfies $N^{\prime}_{\real{}}=N_{\real{}}$.
We will write $N_G:=N'$ when emphasizing that $N'$ is determined by $G$.
\begin{defi}\upshape
Define the {\it age} of an element $g=\frac{1}{r}(a_1,\dots,a_n) \in G$ to be
$$
\age(g)=\frac{1}{r}\sum_{i=1}^n a_i.
$$
Similarly, we define the {\it age} of an element $\bar{g}=\frac{1}{r}(a_1,\dots,a_n) \in N^{\prime}$ to be 
$$
\age(\bar{g})=\frac{1}{r}\sum_{i=1}^n a_i.
$$
\end{defi}
\begin{defi}\upshape
Let $g\in G$ and $I_G$ be the unit of $G$. Then, the rank
$$
\rank(g- I_G)
$$
is called the {\it height} of $g$ and denoted by $\Ht(g)$.
\end{defi}
\begin{prop}[\cite{2}, Prop.5.2.]\label{age-prop}
Let $g\in G$. The following formula holds.
\[\Ht(g)=\Ht(g^{-1}) =\age(g)+\age(g^{-1}).\]
\end{prop}

We shall recall the definition of crepant resolutions in matters of toric geometry.
If a fan $\Sigma$ subdivides the fan $\sigma$, then we have a birational map $f: X(N^{\prime},\Sigma) \to X(N^{\prime},\sigma)$, and the following relation holds between the canonical divisors:
$$
K_{X(N^{\prime},\Sigma)}=f^{*}(K_{X(N^{\prime},\sigma)})+\sum_{\tau \in \Sigma(1)}a_{\tau}D_{\tau},
$$
where $D_{\tau}$ is an exceptional divisor corresponding to the one dimensional cone $\tau \in \Sigma(1)$ in $\Sigma$ and $a_{\tau}=\age(A_{\tau})-1$, where $A_{\tau}$ is the primitive element in $\tau$. The rational number $a_{\tau}$ is called the {\it discrepancy} of $D_{\tau}$. 

\begin{rem}\label{rem2.2}\upshape
Let $\Sigma$ be a subdivision of $\sigma$ by using lattice points whose ages are $1$. If the toric variety $X(N^{\prime},\Sigma)$ is smooth, then $X(N^{\prime},\Sigma)$ is a crepant resolution of $\compl n/G$.
\end{rem}

The convex hull $\mathfrak{s}_G\subset N'_{\real{}}$ spanned by ${\bm e}_1,{\bm e}_2,\ldots ,{\bm e}_n$ is called the {\it junior simplex}. An element in the junior simplex is called a {\it junior element}. By Remark \ref{rem2.2}, a crepant resolution $X(N^{\prime},\Sigma)$ can be identified with a basic triangulation of  $\mathfrak{s}_G$ by using points in $N'$. As for this fact, a necessary condition for Gorenstein abelian quotient singularity to admit a crepant resolution via Hilbert basis is known as follows.

\begin{defi}\upshape(\cite[p.11]{6})
Let $\mathrm{Hlb}_{N'}(\sigma)$ be as follows:
$$
\mathrm{Hlb}_{N'}(\sigma)=\left\{ n \in \sigma \cap (N' \backslash \{ 0 \}) \left|
\begin{array}{c}
 \text{$n$ can not be expresed as}\\  
 \text{the sum of two other vectors}\\ 
 \text{belonging to $\sigma \cap (N'\backslash \{ 0 \})$}
\end{array}
\right.\right\}.
$$
The set $\mathrm{Hlb}_{N'}(\sigma)$ is called the {\it Hilbert basis} of $\sigma$ with reference to the lattice $N'$. 
\end{defi}

\begin{thm}\label{firstexistence}
{\rm (\cite[p.30-31]{6})}
Let $\compl r/G$ be a Gorenstein abelian quotient singularity. If $\mathfrak{s}_G$ has a basic triangulation, then 
$$
\mathrm{Hlb}_{N'}(\sigma)=\mathfrak{s}_G \cap N',
$$
i.e., each of the members of the Hilbert basis of $\sigma$ has to be either a junior element or a vertex of $\mathfrak{s}_G$.
\end{thm}

\subsection{Hirzebruch-Jung Continued Fractions}\label{HJCF}
In this section, we shall mention relations between the minimal resolution of an ${\rm A}_n$ singularity and the Hirzebruch-Jung Continued Fraction obtained from the type of quotient singularities. Let $\compl{2}/G$ be a quotient singularity of $\frac{1}{r}(1,a)$-type where $a \in \inte{}$ and $r\in \natu{}$ are coprime. The {\it Hirzebruch-Jung continued fraction} of $\frac{r}{a}$ is defined as follows:
$$
\frac{r}{a}=x_1-\frac{1}{x_2-\frac{1}{x_3-\cdots\frac{1}{x_s}}}=[x_1,\dots,x_s]
$$
where $x_1,\dots,x_s \in \inte{}_{>0}.$

\begin{figure}[H]
  \begin{center}
   \includegraphics[width=150mm]{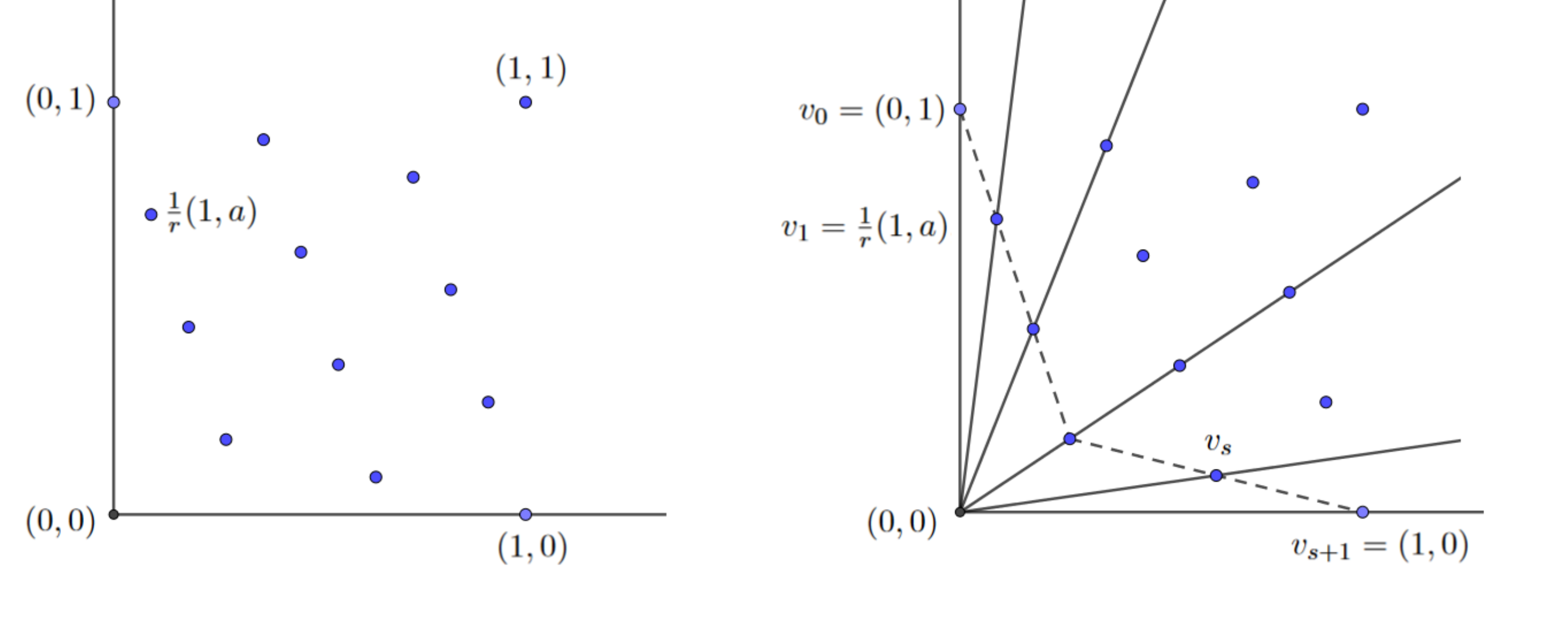}
  \end{center}
  \caption{The minimal resolution of $\compl2/G$ and the Newton Polygon}
  \label{hirzubruchjungex}
\end{figure}

Let $\compl2/G \cong X(N^{\prime},\sigma)$. We set that $\sigma=\real{}_{\geq 0}{\bm e_1}+\real{}_{\geq 0}{\bm e_2}$, $N^{\prime}=\inte2+\inte{}\frac{1}{r}(1,a)$, $v_0={\bm e_2}$ and $v_{s+1}={\bm e_1}$. 
The {\it Newton polygon} $L$ is given as the convex hull of lattice points $\left(N^{\prime}\cap \sigma \right)\setminus \{(0,0)\}$ (See Fig. \ref{hirzubruchjungex}).

\noindent
Let $X(N^{\prime},\Sigma)$ be the minimal resolution of $\compl2/G$. The fan $\Sigma$ is the subdivision of $\sigma$ by the half lines from $(0,0)$ to the primitive elements $v_0={\bm e_2},v_1=\frac{1}{r}(1,a),\dots,v_{s+1}={\bm e_1}$ in $N^{\prime}$. These elements are on the edge of $L$. Moreover, it is known that the following formula holds for the coordinates of these primitive elements and coefficients $x_1,\dots,x_s$ which appear in the Hirzebruch-Jung continued fraction:
$$
v_{i+1}+v_{i-1}=x_iv_i \quad (\ i=1,\dots,s ). 
$$
Therefore, the coordinates of $v_1,\dots,v_s$ can be computed from the Hirzebruch-Jung continued fractions concretely. Every exceptional divisor $E_i$ of the minimal resolution corresponds to the primitive element $v_i$, and its self-intersection number is $-x_i$.

\begin{ex}\upshape
If $a=8$ and $r=11$, then the Hirzebruch-Jung continued fraction is as follows:
$$\frac{11}{8}=2-\frac{1}{2-\frac{1}{3-\frac{1}{2}}}=[2,2,3,2].$$
\noindent
The following list is on the exceptional divisors of the minimal resolution of $X(N^{\prime},\Sigma)$.
\begin{table}[H]
\begin{center}
  \vspace{0pt}\begin{tabular}{|c|c|c|} \hline
    Exceptional Divisors & Primitive Elements in $N^{\prime}$ & Self-Intersection Number  \\  \hline
    $E_1$ & $v_1=\frac{1}{11}\left( 1, 8 \right)$ & $-2$ \\
    $E_2$ & $v_2=\frac{1}{11}\left( 2, 5 \right)$ & $-2$ \\
    $E_3$ & $v_3=\frac{1}{11}\left( 3, 2 \right)$ & $-3$ \\  
    $E_4$ & $v_1=\frac{1}{11}\left( 7, 1 \right)$ & $-2$ \\
\hline
  \end{tabular}
\end{center}
\end{table}
\noindent
Conversely, type of quotient singularities is given by series of coefficients of continued fraction. In the case of $[3,2,2]$, the given quotient singularity is of $\frac{1}{7}(1,3)$-type.

\end{ex}

\subsection{Ashikaga's Continued Fractions}
We shall introduce a generalization of Hirzebruch-Jung continued fractions by Tadashi Ashikaga \cite{1}. This generalized continued fractions summarizes information of the {\it Fujiki-Oka resolution} (See \cite{9, 28}) for {\it semi-isolated} quotient singularities (i.e., cyclic quotient singularities of $\frac{1}{r}(1,a_{2},\ldots,a_{n})$-type. The Fujiki-Oka resolution is a canonical resolution for semi-isolated quotient singularities. We call this continued fraction {\it Ashikaga's continued fractions}.

\begin{defi}\upshape\label{profrac}
Let $n$ be an integer greater than or equal to $1$. Let $\mathbf{a}=(a_1,\dots,a_n) \in \inte n$ and $r \in \natu{}$ which satisfies $0\leq  a_i \leq r-1$ for $1\leq i \leq n$. We call the symbol
$$
\frac{\mathbf{a}}{r}=\frac{(a_1,\dots,a_n)}{r}
$$
an {\it $n$-dimensional proper fraction}.
\end{defi}

\begin{defi}\upshape 
Define the {\it age} of an $n$-dimensional proper fraction $\frac{\mathbf{a}}{r}=\frac{(a_1,\dots,a_n)}{r}$ to be
$$
\age\left( \frac{\mathbf{a}}{r} \right)=\frac{1}{r}\sum_{i=1}^n a_i.
$$
\end{defi}

In the following, the symbol $\ratio{prop}_n$ (resp. $\overline{\ratio{prop}_n}$) means the set of $n$-dimensional proper fractions (resp. the set $\ratio{prop}_n \cup \{\infty\}$). Similarly, $\overline{\inte n}= \inte n \cup \{\infty\}$. Moreover, $\overline{\ratio{prop}_n}[x_2,\dots,x_n]$ (resp. $\overline{\inte n}[x_2,\dots,x_n]$) denotes the set consisting of all noncommutative polynomials with $n-1$ variables over $\overline{\ratio{prop}_n}$ (resp. $\overline{\inte n}$), and $\mathbf{I}=\{2,\dots,n\}$ signifies the index set of the variables $x_2,\dots,x_n$.

Ashikaga's continued fraction consists of a {\it round down polynomial} and a {\it remainder polynomial}, and these polynomials are obtained via {\it round down maps} and {\it remainder maps} for a {\it semi-unimodular proper fraction} (i.e., a proper fraction such that at least one component of $\mathbf{a}$ is $1$). Roughly speaking, these maps are division for just one component of the vector $\mathbf{a}$ by $r$. In the following, we may assume that the first component of a semi-unimodular proper fraction is always $1$ by changing coordinates.\label{ACF}

\begin{defi}\upshape (\cite[Def 3.1.]{1})
Let $\frac{(1,a_2,\dots,a_n)}{r}$ be a semi-unimodular proper fraction. 
\begin{itemize}
\item[$(i)$] For $2\leq i \leq n$, the {\it $i$-th remainder map} $R_i:\overline{\ratio{prop}_n} \to \overline{\ratio{prop}_n}$ is define by
$$
R_i\left(\frac{(1, a_2,\dots,a_n)}{r}\right)=
\left\{ 
\begin{array}{cc}
 \frac{\left(\overline{1}^{a_i}, \overline{a_2}^{a_i},\ \dots,\ \overline{a_{i-1}}^{a_i},\ \overline{-r}^{a_i},\ \overline{a_{i+1}}^{a_i},\ \dots, \overline{a_n}^{a_i}\right)}{a_i}
         & {\rm if} \  a_i \neq 0\\  
 \infty &  {\rm if}\  a_i=0
\end{array}
\right.
$$
and $R_i(\infty)=\infty$ where $\overline{a_j}^{a_i}$ is an integer satisfying $0\leq \overline{a_j}^{a_i} < a_i$ and $\overline{a_j}^{a_i} \equiv a_j$ modulo $a_i$.
\item[$(ii)$] For $2\leq i \leq n$, the {\it $i$-th round down map} $Z_i:\overline{\ratio{prop}_n} \to \overline{\inte n}$ is defined by
$$
Z_i\left(\frac{(1,a_2,\dots,a_n)}{r}\right)=
\left\{ 
\begin{array}{cc}
 \left( \lfloor \frac{1}{a_i} \rfloor, \lfloor \frac{a_2}{a_i} \rfloor, \dots, \lfloor \frac{a_{i-1}}{a_i}\rfloor,  \lfloor \frac{-r}{a_i}\rfloor,  \lfloor \frac{a_{i+1}}{a_i}\rfloor, \dots,  \lfloor \frac{a_n}{a_i}\rfloor \right)& {\rm if} \  a_i \neq 0\\  
 \infty &  {\rm if}\  a_i=0
\end{array}
\right.
$$
and $Z_i(\infty)=\infty$ where $ \lfloor x \rfloor$ is the greatest integer not exceeding $x$.
\end{itemize}
\end{defi}

\begin{ex}\upshape
If $v=\frac{(1,2,5,7)}{8}$, then 
\begin{eqnarray*}
Z_2(v)&=&(0,-4,2,3), \\
Z_3(v)&=&(0,0,-2,1), \\
R_2(v)&=&\frac{(1,0,1,1)}{2}\ \text{and}\\
R_3(v)&=&\frac{(1,2,2,2)}{5}.
\end{eqnarray*}
\end{ex}

\begin{defi}\upshape\cite[Def 3.2.]{1}\label{DOACF}
Let $\frac{\mathbf{a}}{r}$ be an $n$-dimensional semi-unimodular proper fraction. 
\begin{itemize}
  \item[(i)] The {\it remainder polynomial} $\mathcal{R}_*\left(\frac{\mathbf{a}}{r}\right) \in \overline{\ratio{prop}_n}[x_2,\dots,x_n]$ is defined by
  $$
  \mathcal{R}_*\left(\frac{\mathbf{a}}{r}\right)=\frac{\mathbf{a}}{r}+
                                      \sum_{(i_1,i_2,\dots,i_l)\in \mathbf{I}^l,\: l\geq 1 }(R_{i_l}\cdots R_{i_2}R_{i_1})\left(\frac{\mathbf{a}}{r}\right)\cdot x_{i_1}x_{i_2}\cdots x_{i_l}
  $$
  where we exclude terms with coefficients $\infty$ or $\frac{(0,0,\dots,0)}{1}$.
  \item[(ii)] The {\it round down polynomial} $Z_* \in \overline{\inte n}[x_2,\dots,x_n]$ is defined by
  $$
  \mathcal{Z}_*\left(\frac{\mathbf{a}}{r}\right)=\sum_{j=2}^n Z_j\left(\frac{\mathbf{a}}{r}\right)x_j+
                                       \sum_{j=2}^n  \sum_{(i_1,i_2,\dots,i_l)\in \mathbf{I}^l,\: l\geq 1 }(Z_jR_{i_l}\cdots     R_{i_2}R_{i_1})\left(\frac{\mathbf{a}}{r}\right)\cdot x_{i_1}x_{i_2}\cdots x_{i_l}x_j.
  $$
 \end{itemize}
\end{defi}

\begin{rem}\upshape
In the case $n=2$, the series of the coefficients of $\mathcal{Z}_*\left(\frac{\mathbf{a}}{r}\right)$ coincides with the series of the coefficients of Hirzebruch-Jung continued fraction.
\end{rem}

\begin{ex}\upshape
Let $v=\frac{(1,2,8)}{11}$, then the remainder polynomial is
\begin{eqnarray*}
\mathcal{R}_*\left(\frac{(1,2,8)}{11}\right)=
             \frac{1}{11}(1,2,8)&+&\frac{1}{2}(1,1,0)x_2+\frac{1}{8}(1,2,5)x_3 \\
                                &+&\frac{1}{2}(1,0,1)x_3x_2+\frac{1}{5}(1,2,2)x_3x_3 \\
                                &+&\frac{1}{2}(1,1,0)x_3x_3x_2+\frac{1}{2}(1,0,1)x_3x_3x_3.
\end{eqnarray*}
The round down polynomial is
\begin{eqnarray*}
\mathcal{Z}_*\left(\frac{(1,2,8)}{11}\right)&=&
             (0,-6,4)x_2+(0,0,-2)x_3 \\
                                &+&(1,-4,2)x_3x_2+(0,0,-2)x_3x_3 \\
                                &+&(0,-3,1)x_3x_3x_2+(0,1,-3)x_3x_3x_3.
\end{eqnarray*}
\end{ex}

Remainder polynomials consist of datum of blow-up centers of a Fujiki-Oka resolution. In the following, we shall summarize Fujiki-Oka resolutions. For the details of the Definition \ref{unimodular} and the Lemma \ref{Okalem}, see the articles written by T. Ashikaga \cite{1} and M. Oka \cite{28} respectively.

\begin{defi}\upshape\label{unimodular}
Let $P_1,\ldots,P_n$ be primitive elements in $N^{\prime}$. If an $n$-dimensional cone $\tau=\real{}_{\geq 0}P_1+\dots+\real{}_{\geq 0}P_n$ in $N^{\prime}_{\real{}}$ has a smooth facet $\real{}_{\geq 0}P_1+\dots+\real{}_{\geq 0} P_{n-1}$, then we call the cone {\it semi-unimodular} over the vertex $P_n$.
\end{defi}

If a cone $\tau$ is semi-unimodular over all vertices $P_1,\ldots,P_n \in N^{\prime}$, then the toric variety $X(N^{\prime},\tau)$ has an isolated singularity or no singularities. If the toric variety $X(N^{\prime},\sigma)$ has a quotient singularity of $\frac{1}{r}(a_1,\dots,a_n)$-type satisfying $\mathrm{GCD}(r,a_i)=1$, then $\sigma$ is semi-unimodular over ${\bm e}_i$, and $X(N^{\prime},\sigma)$ has a semi-isolated singularity.

\begin{lem}\label{Okalem}
Let an $n$-dimensional cone $\tau=\real{}_{\geq 0}P_1+\dots+\real{}_{\geq 0}P_n \subset N^{\prime}_{\real{}}$ be semi-unimodular over $P_1$ and $r=\det(P_1,P_2,\cdots, P_n)$. If $C\in \inte{n}$ is a primitive element such that the $n$-dimensional cone $\real{}_{\geq 0}C+\real{}_{\geq 0}P_2+\dots+\real{}_{\geq 0}P_n$ is smooth, then there exist integers $0 \leq a_2,\dots,a_n \leq r-1$ such that 
 $$
 C=\frac{P_1+\sum_{i=2}^{n}a_iP_i}{r}.
 $$
\end{lem}
This element $C\in N^{\prime}$ is called a {\it Oka center} of $\tau$ over $P_1$. The Oka center exists uniquely for a cone which is semi-unimodular over an element.
\begin{lem}{\rm (\cite[Lemma 3.]{9})}\label{Fujilem}
Suppose $(X,[0])$ is a cyclic quotient singularity of $\frac{1}{r}(a_1,a_2,\dots,a_n)$-type where $\mathrm{GCD}(r,a_1,\ldots,a_n)=1$, and $a_1 a_2 \cdots a_n \neq 0$. Then there exist a variety $\widetilde{X}$, a finite affine open covering $\mathcal{U} = \{ U_1, \ldots,U_l \}$ of $\widetilde{X}$ for an integer $1\leq l \leq n$, and a proper birational morphism $f: \widetilde{X}\to X$ such that $U_i$ is the quotient singularity of $\frac{1}{r_i}(a_{i1},a_{i2},\dots,a_{in})$-type for each $i$, where the integers $r_i$ and $a_{ij}\ (1\leq j\leq n)$ are determined by the following formula:
$$
\left\{ 
\begin{array}{l}
 r_i = a_i / d\ \text{where}\ d=\mathrm{GCD}(a_1, \ldots, a_n),\\  
 a_{ij} \equiv r_j\ \text{modulo}\ r_i\ \text{and}\ 0\leq a_{ij}<r_i\ ( j\neq i),\\
 a_{ij} + r \equiv 0\ \text{modulo}\ r_i\ ( j= i).
\end{array}
\right.
$$
\end{lem}

The proper fraction $\frac{\mathbf{a}}{r}=\frac{(1,a_2,\ldots,a_n)}{r}$ obtained from Lemma \ref{Okalem} and Lemma \ref{Fujilem} is called the {\it proper fraction of $\tau$ over $P_1$}.

\begin{lem}{\rm(T. Ashikaga, \cite{1})}\label{Ashilem}
If a cone $\tau=\real{}_{\geq 0}P_1+\dots+\real{}_{\geq 0}P_n \subset N^{\prime}_{\real{}}$ contains a primitive element $C\in N^{\prime} $ in Lemma \ref{Okalem}, then the toric variety $X(N^{\prime},\tau)$ has a quotient singularity of $\frac{1}{r}(1,a_2,\dots,a_n)$-type.
\end{lem}

By Lemma \ref{Okalem}, Lemma \ref{Fujilem} and Lemma \ref{Ashilem}, any semi-isolated quotient singularity is resolved by blow-ups with the Oka center repeatedly, and these toric resolutions are called {\it Fujiki-Oka resolutions}. Each coefficient which appears in a remainder polynomial coincides with the type of quotient singularities appearing in each step of the Fujiki-Oka resolution. 

\begin{lem}{\rm (T. Ashikaga, \cite{1})}\label{Ashilem2}
Let $\tau=\real{}_{\geq 0}P_1+\dots+\real{}_{\geq 0}P_n \subset N^{\prime}_{\real{}}$ be a semi-unimodular cone over $P_1$ and $C$ be the Oka center of $\tau$. Then the cone 
$$
\tau_i=\real{}_{\geq 0}P_1+\dots+\real{}_{\geq 0}P_{i-1}+\real{}_{\geq 0}c_i+\real{}_{\geq 0}P_{i+1}+\dots+\real{}_{\geq 0}P_n
$$
is semi-unimodular over $P_1$ and its Oka center is 
$$
c_i=\frac{\sum_{j\neq i,n}\overline{a_j}^{a_i}P_j+\overline{-d}^{a_i}P_i}{a_i}.
$$
\end{lem}
By Lemma \ref{Ashilem2}, remainder polynomials can be understood as a overview of Fujiki-Oka resolution of a semi-isolated quotient singularity.

\begin{ex}\upshape
Let $X(N^{\prime},\sigma)$ have a quotient singularity of $\frac{1}{11}(1,2,8)$-type, i.e., $N^{\prime}=\inte3+\inte{}\frac{1}{11}(1,2,8)$ and $\sigma=\real{}_{\geq 0}{\bm e}_1+\real{}_{\geq 0}{\bm e}_2+\real{}_{\geq 0}{\bm e}_3$. Then, the cone $\sigma$ is semi-unimodular over ${\bm e}_1$, and the Oka center is $c=\frac{1}{11}(1,2,8)$, and the remainder polynomial of the proper fraction $\frac{(1,2,8)}{11}$ is as follows:
\begin{eqnarray}
\mathcal{R}_*\left(\frac{(1,2,8)}{11}\right)=
             \frac{1}{11}(1,2,8)&+&\frac{1}{2}(1,1,0)x_2+\frac{1}{8}(1,2,5)x_3 \nonumber \\
                                &+&\frac{1}{2}(1,0,1)x_3x_2+\frac{1}{5}(1,2,2)x_3x_3 \nonumber \\
                                &+&\frac{1}{2}(1,1,0)x_3x_3x_2+\frac{1}{2}(1,0,1)x_3x_3x_3.\nonumber
\end{eqnarray}
This expanding of Ashikaga's continued fraction indicates that the toric variety after the blow-up with the Oka center $\frac{1}{11}(1,2,8)$ has two semi-isolated quotient singularities of $\frac{1}{2}(1,1,0)$-type and $\frac{1}{8}(1,2,5)$-type. For these quotient singularities, the corresponding cones which appear in $\sigma$ after the subdivision by $\frac{1}{11}(1,2,8)\in N^{\prime}$ are $\sigma_2=\real{}_{\geq 0}{\bm e}_1+\real{}_{\geq 0}c+\real{}_{\geq 0}{\bm e}_3$ and $\sigma_3=\real{}_{\geq 0}{\bm e}_1+\real{}_{\geq 0}{\bm e}_2+\real{}_{\geq 0}c$ respectively. $\frac{1}{2}(1,1,0)$ and $\frac{1}{8}(1,2,5)$ are the Oka center of semi-unimodular cones $\sigma_2,\ \sigma_3$ over ${\bm e}_1$ respectively. Therefore, we can take blow-ups with the Oka centers again. The blow-up with Oka centers of $X(N^{\prime},\sigma_3)$ consists a smooth toric variety and quotient singularities of $\frac{1}{2}(1,0,1)$-type and $\frac{1}{5}(1,2,2)$-type respectively. By repeating blow-ups with Oka centers, we have the smooth toric variety (See Fig. \ref{128}). 

\begin{figure}[H]
  \begin{center}
   \includegraphics[width=150mm]{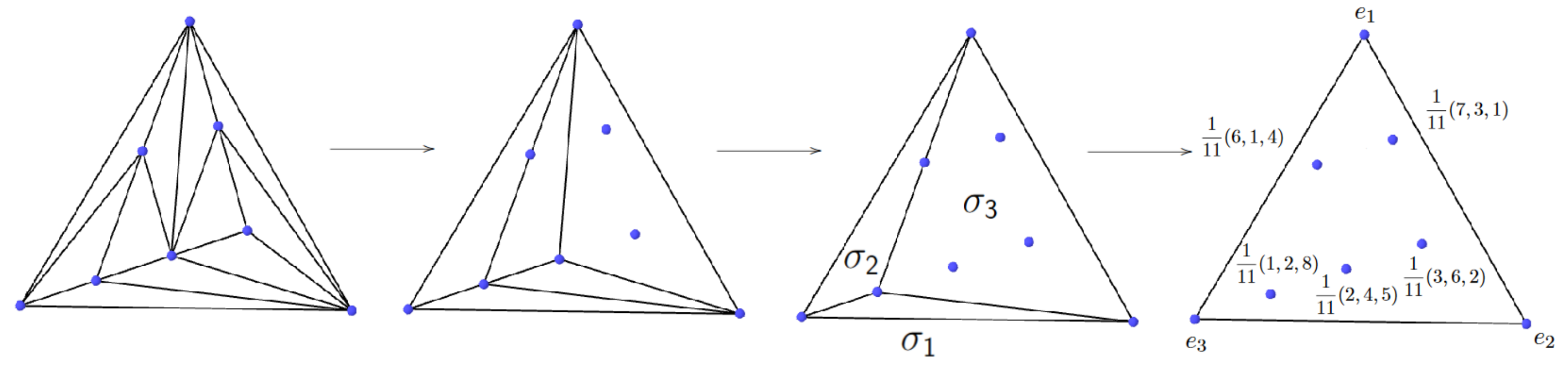}
  \end{center}
  \caption{ The basic triangulation of $\mathfrak{s}_G$ by Fujiki-Oka resolution}
  \label{128}
\end{figure}

\end{ex}
Let $X(N^{\prime},\sigma)$ have a semi-isolated quotient singularity. The cone $\sigma \subset N^{\prime}_{\real{}}$ can be semi-unimodular over ${\bm e}_1$ by exchanging basis of $N^{\prime}$. For a semi-unimodular cone $\sigma$ over ${\bm e}_1$, we call the terminal smooth fan obtained from its Fujiki-Oka resolution the {\it continued fraction fan}, and that fan is denoted as $\CFF_{{\bm e}_1} (\sigma)$ or, more simply, $\CFF (\sigma)$. Clearly, there exists at least one $\CFF (\sigma)$ for a semi-unimodular cone $\sigma$.  

\section{Crepant resolutions via Ashikaga's continued fractions}
In this section, we shall show a sufficient condition of existence of a crepant resolution for Gorenstein abelian quotient singularities by using Ashikaga's continuous fractions as the main result of this paper. 

\subsection{Sufficient condition of existence of a crepant resolution}
The purpose of this subsection is to show a sufficient condition of existence of a crepant resolution for semi-isolated cyclic quotient singularities. In particular, all isolated cyclic quotient singularities are included in this case.

\begin{thm}\label{thm1}
For a cyclic quotient singularity of $\frac{1}{r}(1,a_2,\dots,a_n)$-type, the Fujiki-Oka resolution is crepant if and only if the ages of all the coefficients of the corresponding remainder polynomial $\mathcal{R}_*\left(\frac{(1,a_2,\dots,a_n)}{r}\right)$ are $1$.
\end{thm}
Let $G=\left\langle\frac{1}{r}(1,a_2,\dots,a_n)\right\rangle=\langle g\rangle$, we will denote by $G_i$ the cyclic group which is generated by $g_i$, where $g_i$ is determined by the image of the $i$-th remainder map 
$$
 \mathcal{R}_i\left(\frac{(1,a_2,\dots,a_n)}{r}\right) = \frac{(1,\overline{a_2}^{a_i},\dots,\overline{a_{i-1}}^{a_i},\overline{-r}^{a_i},\overline{a_{i+1}}^{a_i},\dots,\overline{a_n}^{a_i})}{a_i},
$$
i.e., $\compl{n}/G_i$ is the cyclic quotient singularity of $\frac{1}{a_i}(1,\overline{a_2}^{a_i},\dots,\overline{a_{i-1}}^{a_i},\overline{-r}^{a_i},\overline{a_{i+1}}^{a_i},\dots,\overline{a_n}^{a_i})$-type.
We first show that the generator $g=\frac{1}{r}(1,a_2,\dots,a_n)$ satisfies $\age(g)=1$ if $\compl{n}/G$ has a crepant resolution.

\begin{prop}\label{Prop1}
Assume that $1+a_2+\dots+a_n \geq 2r$ for $G=\left\langle\frac{1}{r}(1,a_2,\dots,a_n)\right\rangle$. Then, $\compl n/G$ has no toric crepant resolutions.
\end{prop}
\proofname
Assume that $\compl n/G \cong X(\sigma, N^{\prime})$ has a toric crepant resolution $X(\Sigma, N^{\prime})$.
Since $G$ has the generator of which the first component is $\frac{1}{r}$, we see that there are no lattice points on $\mathfrak{s}_G \cap \tau_1$ where $\tau_1$ is  $n-1$ dimensional cone with vertices ${\bm e}_2,\dots,{\bm e}_n$.
Therefore, there is a lattice point ${\bm q} \in N'$ such that $\age({\bm q})=1$ and ${\rm Cone}({\bm q}, {\bm e}_2, {\bm e}_3,\dots ,{\bm e}_n) \in \Sigma$.
Since $G$ is generated by $\frac{1}{r}(1,a_2,\dots,a_n)$ with  $1+a_2+\dots+a_n \geq 2r$
, we can write ${\bm q}=\frac{1}{r}(i,\overline{a_2i}^r,\overline{a_3i}^r,\dots,\overline{a_ni}^r)$ where $i \neq 1$.
Thus, we have $\{{\bm q}, {\bm e}_2, {\bm e}_3,\dots ,{\bm e}_n\}$ as $\inte{}$-basis of $N'$, so there exist integers $k_1,\dots,k_n \in \inte{}$ such that

$$
{\bm p}=\frac{1}{r}(1,a_2,\dots,a_n)=k_1{\bm q}+\sum_{j=2}^nk_j{\bm e}_j.
$$

We now turn to the first component. This formula gives  $\frac{1}{r}=\frac{k_1i}{r}$, but it contradicts $i\neq1$.
Therefore, if $1+a_2+\dots+a_n \geq 2r$, then $\compl n/G$ has no crepant resolutions. \qed\\
\thefootnote

\begin{prop}\label{Prop2}
Let $G=\left\langle\frac{1}{r}(1,a_2,\dots,a_n)\right\rangle$ with $1+a_2+\dots+a_n=r$. If $\compl n/G_i$ have a crepant resolution for all $i=2,\dots,n$ , then $\compl n/G$ have a crepant resolution. 
\end{prop}
\proofname \ 
Let $X(\Sigma_i, N_i)$ be a toric crepant resolution for $\compl n/G_i$ where $N_i=\inte n + g_i\inte{}$ with canonical basis $\overline{{\bm e}_1}\dots,\overline{{\bm e}_n}$. For simplicity of notation, we write $N_{i\mathbb{R}}$ insteads of $N_i\otimes_{\inte{}} \mathbb{R}$.
Fix a smooth cone in $\Sigma_i$, and write this cone $\sigma={\rm Cone}(\overline{{\bm v}_1},\dots,\overline{{\bm v}_n})$. For ${\bm x}=(x_1, \dots, x_n) \in N_{i\mathbb{R}}$, the map $\phi_i:N_{i\mathbb{R}} \hookrightarrow N'_{\mathbb{R}}$ is defined as follows:
$$
\phi_i({\bm x})=\left(x_1+\frac{1}{r}x_i, x_2+\frac{a_2}{r}x_i,\dots,  x_{i-1}+\frac{a_{i-1}}{r}x_i, \frac{a_i}{r}x_i, x_{i+1}+\frac{a_{i+1}}{r}x_i, \dots, x_n+\frac{a_n}{r}x_i \right).
$$
The proof will be divided into two steps. The first step is to check $\phi_i({\bm x})$ satisfies $\age(\phi_i({\bm x}))=1$ for a point ${\bm x} \in N_{i\mathbb{R}}$ with $\age({\bm x})=1$, the second step is to prove $\phi_i(\sigma) \subset N'_{\mathbb{R}}$ is also smooth.
\begin{itemize}
\item[(i)] If ${\bm x}=(x_1,\dots,x_n) \in N_i$ satisfies $\age({\bm x})=1$, then ${\rm age}(\phi_i(x))=x_1+x_2+\dots+\hat{x_i}+\dots+x_n+\frac{1}{r}(1+a_2+\dots+a_n)x_i$. By assumption, we have $x_1+\dots+x_n=1$ and $1+a_2+\dots+a_n=r$. These formula give ${\rm age}(\phi_i({\bm x }))=x_1+\dots+x_n=1$．

\item[(ii)] Since $\sigma$ is smooth on $N_i$, $\{\overline{{\bm v}_1},\dots,\overline{{\bm v}_n}\}$ is $\inte{}$-basis of $N_i$, namely  genarate the canonical basis $\overline{{\bm e}_1},\dots,\overline{{\bm e}_n}$ of $N_i$ and $\overline{{\bm q}}=\frac{1}{a_i}(1,\overline{a_2}^{a_i},\dots,\overline{-r}^{a_i},\dots, \overline{a_n}^{a_i})$. Let ${\bm v}_j$ denote $\phi_i(\overline{{\bm v}_j})$, then we have $\phi_i(\sigma)={\rm Cone}({\bm v}_1,\dots,{\bm v}_n)$. It is easy to see that $\phi_i({\bm e}_j)$ and $\phi_i(\overline{{\bm q}})$ are generated by $V=\{{\bm v}_1,\dots,{\bm v}_n\}$, where $j \in \{1,\dots,n\} \backslash\{i\}$. 
To show $V$ is $\inte{}$-basis of $N'$, it is sufficient to prove that ${\bm e}_i$ is generated by $V$. Let us denote by $Q_z$ the quotient of $z$ devided by $a_i$. We have
$$
{\bm q}=\phi_i(\overline{{\bm q}})=\frac{1}{r}(-Q_{-r},Q_{a_2}r-a_2Q_{-r},\dots,(-r)-a_iQ_{-r},\dots,Q_{a_n}r-a_nQ_{-r}),
$$
and we get the formula
$$
{\bm q}+Q_{-r}{\bm p}=\frac{1}{r}(0,Q_{a_2}r,\dots,Q_{a_{i-1}}r, -r, Q_{a_{i+1}}r,\dots,Q_{a_n}r).
$$
Therefore, the following equation holds
$$
{\bm e}_i={\bm q}+Q_{-r}{\bm p}-\sum_{j \in \{1,\dots,n\}\backslash\{i\}} Q_{a_j}{\bm e}_j.
$$
This implies that $\{{\bm v}_1,\dots,{\bm v}_n\}$ generates ${\bm e}_i$. Thus, $\{{\bm v}_1,\dots,{\bm v}_n\}$ is $\inte{}$-basis of $N'$.
\end{itemize}

From  {\rm (i)} and {\rm (ii)}, we see that if $\compl n/G_i$ has a  crepant resolution $X(\Sigma_i,N_i)$, then we have a fan on $N'$ corresponding to a crepant resolution for $\compl n/G$ by taking the union of all $\phi_i(\Sigma_i)$.
\qed\\

\textbf{The proof of Theorem\ref{thm1}:} 
 Assume that the ages of all coefficients of the remainder polynomial of $\frac{(1,a_2,\dots,a_n)}{r}$ equal to $1$. By Proposition\ref{Prop2}, whether $\compl n/G$ has a crepant resolution depends on whether $\compl n/G_i$ has a crepant resolution for all $i$. It is obvious that the order of $G_i$ is less than the order of $G$. The repeated application of the remainder map enables us to get $G_{i_1i_2\cdots i_j}=\frac{(1,c_2,c_3,\dots,c_n)}{k}$ with $c_j \in \{0,1\}$ for all $j$.
Since $\age\left(\frac{1}{k}(1,c_2,c_3,\dots,c_n)\right)=1$, the Fujiki-Oka resolution of $\compl n/G_{i_1i_2\cdots i_j}$ is crepant. By the proof of Proposition\ref{Prop2}, the Fujiki-Oka resolution of $\compl n/G$ is crepant. Conversely, if the Fujiki-Oka resolution of $\compl n/G$ is crepant, then $\age\left(\frac{1}{r}(1,a_2,a_3,\dots,a_n)\right)=1$ and $\age(g_i)=1$ for $i=2,\dots, n$. Therefore, the ages of all the coefficients of the remainder polynomial are $1$, which completes the proof.\\

The Gorenstein property of $\compl n/G_i$ comes from the property of the cyclic group $G$. We have the following Lemma.

\begin{lem}\label{Lem1}
Assume that $1+a_2+a_3+\dots+a_n=r$ for $G=\left\langle \frac{1}{r}(1,a_2,\dots,a_n) \right\rangle$. Then $\age\left(\mathcal{R}_i\left(\frac{(1,a_2,\dots,a_n)}{r}\right)\right)$ is an integer.
\end{lem}
\proofname
It is enough to prove in the case of $i=2$. We have the equation 
$\mathcal{R}_2\left(\frac{(1,a_2,\dots,a_n)}{r}\right)=\frac{(1,\overline{-r}^{a_2},\overline{a_3}^{a_2},\dots,\overline{a_n}^{a_2})}{a_2}.$
We claim that $1+\overline{-r}^{a_2}+\overline{a_3}^{a_2}+\dots+\overline{a_n}^{a_2}$ is divided by $a_2$.
It is sufficient to show that $1+(-r)+a_3+a_4+\dots+a_n$ is divided by $a_2$. By the assumption $1+a_2+a_3+\dots+a_n=r$, we have $1+(-r)+a_3+a_4+\dots+a_n=a_2$. Therefore, $\age\left(\mathcal{R}_2\left(\frac{(1,a_2,\dots,a_n)}{r}\right)\right)$ is an integer.
\qed\\

Lemma \ref{Lem1} and Theorem \ref{thm1} lead to the following corollary.

\begin{cor}\label{maincor}
For all three dimensional semi-isolated Gorenstein quotient singularities, the Fujiki-Oka resolutions are crepant.
\end{cor}

\proofname Let $G=\left\langle \frac{1}{r}(1,a,b)\right \rangle$ where $1+a+b=r$. we have $\mathcal{R}_2\left(\frac{(1,a,b)}{r}\right)=\frac{(1,\overline{-r}^a,\overline{b}^a)}{a}$, and the age of $\mathcal{R}_2\left(\frac{(1,a,b)}{r}\right)$ is an integer by Lemma \ref{Lem1}.
Clearly, $1+\overline{-r}^a+\overline{b}^a< 2a $. So, the age of $\mathcal{R}_2\left(\frac{(1,a,b)}{r}\right)$ equals to $1$.
Thus, the ages of all coefficients of $\mathcal{R}_*\left(\frac{(1,a,b)}{r}\right)$ equal to $1$. By Theorem \ref{thm1}, the Fujiki-Oka resolution $X(N^{\prime},\CFF(\sigma))$ is crepant.
\qed\\

\subsection{First Existence Criterion via Continued Fractions}

We will give the continued fraction version of Theorem \ref{firstexistence}. 

\begin{defi}\upshape
The term with the variable $x_i\cdots x_i$ in a remainder polynomial is called {\it iterated} where $1 \leq i \leq n$, and 
the lattice point in $N'$ corresponding to the coefficient of iterated terms is also called to be {\it iterated}.
\end{defi}

Every iterated point can be written as $\phi_i^{-1}(\frac{\mathbf{a}}{r}) \in N'$ for the coefficient $\frac{\mathbf{a}}{r}$ of an iterated term.

We shall consider a relationship between iterated points and Hilbert basis, and apply the relationship to Theorem \ref{firstexistence}. 

In the following, for a cyclic group $\left\langle\frac{1}{r}(1,a_2,\dots,a_n)\right\rangle \subseteq \slmc n$ satisfying $1+a_2+\cdots+a_n=r$, the symbol $A_i$ denotes the cyclic subgroup $\left\langle\frac{1}{r}(1,a_i)\right\rangle \subset\glmc2$ for $a_i \neq 0$, and ${\bm v}_{i_1},\dots,{\bm v}_{i_s}$ denote the lattice points in $N_{A_i}=\inte2+\frac{1}{r}(1,a_i)\inte{}$ such that 
$$
{\bm v}_{i_{j-1}}+{\bm v}_{i_{j+1}}=\alpha_{i_{j}} {\bm v}_{i_j}\ \text{for}\ j=1,\dots,s 
$$
where the integers $\alpha_{i_1},\dots, \alpha_{i_s}$ are the entries of the Hirzbruch-Jung continued fraction $\frac{r}{a_i}=[\alpha_{i_1},\dots,\alpha_{i_s}]$ and ${\bm v}_{i_0}=(0,1), {\bm v}_{i_{s+1}}=(1,0)$.
The lattice point ${\bm v}_{i_j}$ can be written as ${\bm v}_{i_j}=\frac{1}{r}(k_{i_j},\overline{a_i\cdot k_{i_j}}^r)$ for some positive integers $k_{i_j}$.

\begin{defi}\upshape
Let $r,\ a_i$ and $k_{i_j}$ be as above. We define an {\it $i$-th minimal point} ${\bm u}_{ij}\in N^{\prime}$ as follows:
$${\bm u}_{ij}=\frac{1}{r}(k_{i_j},\overline{a_2\cdot k_{i_j}}^r,\dots,\overline{a_n\cdot k_{i_j}}^r).$$ 
\end{defi}

We note that ${\bm v}_{i_0}, \dots,{\bm v}_{i_{s+1}}$ are elements in $\mathrm{Hlb}_{N_{A_i}}(\sigma_{A_i})$, where $\sigma_{A_i}=\mathrm{Cone}((1,0),(0,1)) \subset N_{A_i}\otimes \real{}$. One of the good properties of minimal points is that they are in Hilbert basis as shown in the next lemma.

\begin{lem}
All minimal points are in $\mathrm{Hlb}_{N'}(\sigma)$.
\end{lem}
\proofname 
Let ${\bm u}=(u_1,\dots,u_n)\in N'$ be an $i$-th minimal point and ${\bm v}=(u_1,u_i)\in N_{A_i}$ be the element corresponding to ${\bm u}$.
If ${\bm u} \notin \mathrm{Hlb}_{N'}(\sigma)$ , then there 
exists $X=(x_1,\dots,x_n)$ and $Y=(y_1,\dots,y_n)$ in $N'$ such that ${\bm u}=X+Y$.
By focusing on the first and $i$-th components, the following equations hold:
\begin{eqnarray*}
u_1&=&x_1+y_1, \\
u_i&=&x_i+y_i.
\end{eqnarray*}
Let $X_i=(x_1,x_i),\ Y_i=(y_1,y_i) \in N_{A_i}$, then we have ${\bm v}=X_i+Y_i$ by the above formula. This contradicts the fact that ${\bm v} \in \mathrm{Hlb}_{N_{A_i}}(\sigma_{A_i})$.
Therefore, we get ${\bm u} \in \mathrm{Hlb}_{N'}(\sigma)$. \qed\\

An iterated point is either a minimal point or a sum of canonical basis and a minimal point. An iterated point is minimal if and only if it satisfies the conditions of the proper fractions. See Definition \ref{profrac}.

\begin{prop}\label{Prop3}
Let $\compl n/G$ be a quotient singularity of $\frac{1}{r}(1,a_2,\dots,a_n)$-type satisfying $1+a_2+\cdots+a_n=r$. If the remainder polynomial $\mathcal{R}_*\left(\frac{(1,a_2,\dots,a_n)}{r}\right)$ contains an iterated term of which the age of the coefficient is equal to or larger than $2$, then $\compl n/G$ has no toric crepant resolutions.
\end{prop}

The problem with Proposition \ref{Prop3} is that if $G$ has some representation $\frac{1}{r}(1,a_2,\dots,a_n)$, $\frac{1}{r}(b_1,1,b_3,\dots,b_n)$ and so on, their remainder polynomials are different from each other, so if necessary, we have to calculate iterated points for all representation. Moreover, in higher dimension, there are many groups which fulfil Theorem \ref{firstexistence} and possess no crepant resolutions. For example, $G=\left\langle\frac{1}{39}(1,5,8,25)\right\rangle$.

\subsection{Examples}
We shall see some examples of continued fractional fans in dimension four. In this subsection, a proper fraction $\frac{(a_1,\ldots,a_n)}{r}$ and the symbol of types of singularities $\frac{1}{r}(a_1,\ldots,a_n)$ are identified with reach other if necessary.

\begin{ex}\upshape
For the isolated quotient singularity of $\frac{1}{15}(1,2,4,8)$-type, the remainder polynomial is;

\begin{eqnarray*}
&\mathcal{R}_*\left(\frac{(1,2,4,8)}{15}\right)&=\frac{(1,2,4,8)}{15}+\frac{(1,1,0,0)}{2}x_2+\frac{(1,2,1,0)}{4}x_3+\frac{(1,2,4,1)}{8}x_4 \\
&\quad&+\frac{(1,0,1,0)}{2}x_3x_2+\frac{(1,0,0,1)}{2}x_4x_2+\frac{(1,2,0,1)}{4}x_4x_3+\frac{(1,0,0,1)}{2}x_4x_3x_2,
\end{eqnarray*}

where we exclude the terms with the coefficient $\frac{(0,0,0,0)}{1}$.
Since the ages of all the proper fractions in the remainder polynomial are $1$, the continued fractional fan $\CFF_{{\bm e}_1}(\sigma)$ gives a crepant resolution $X(N^{\prime},\CFF(\sigma))$ in this case.

Let us see the first subdivision of $\sigma$ by using the point $\frac{1}{15}(1,2,4,8)\in N'$. After the subdivision, the junior simplex $\mathfrak{s}_G$ is divided into the four tetrahedrons as Fig. \ref{1248-3}. The colored tetrahedron in Fig. \ref{1248-3} corresponds to the smooth, and the other three tetrahedrons correspond to the cyclic quotient singularities of $\frac{1}{4}(1,2,1,0)$-type, $\frac{1}{8}(1,2,4,1)$-type and $\frac{1}{2}(1,1,0,0)$-type in order from the left side in Fig. \ref{1248-2}.

\begin{figure}[htbp]
 \begin{minipage}{0.45\hsize}
  \begin{center}
   \includegraphics[width=60mm]{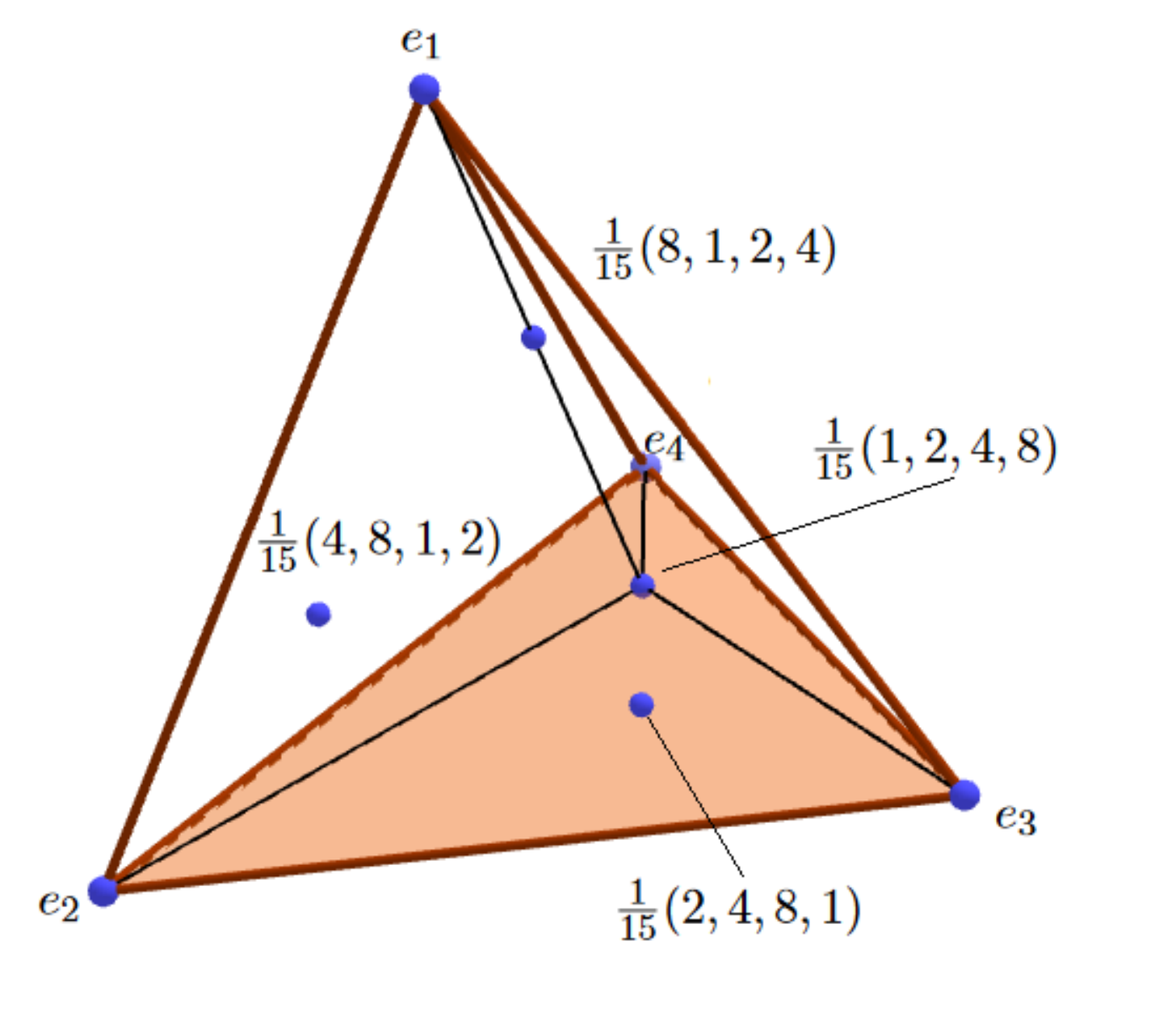}
  \end{center}
  \caption{The junior simplex $\mathfrak{s}_G$ in $N'_{\real{}}$}
   \label{1248-3}
 
 \end{minipage}
 \begin{minipage}{0.55\hsize}
  \begin{center}
   \includegraphics[width=73mm]{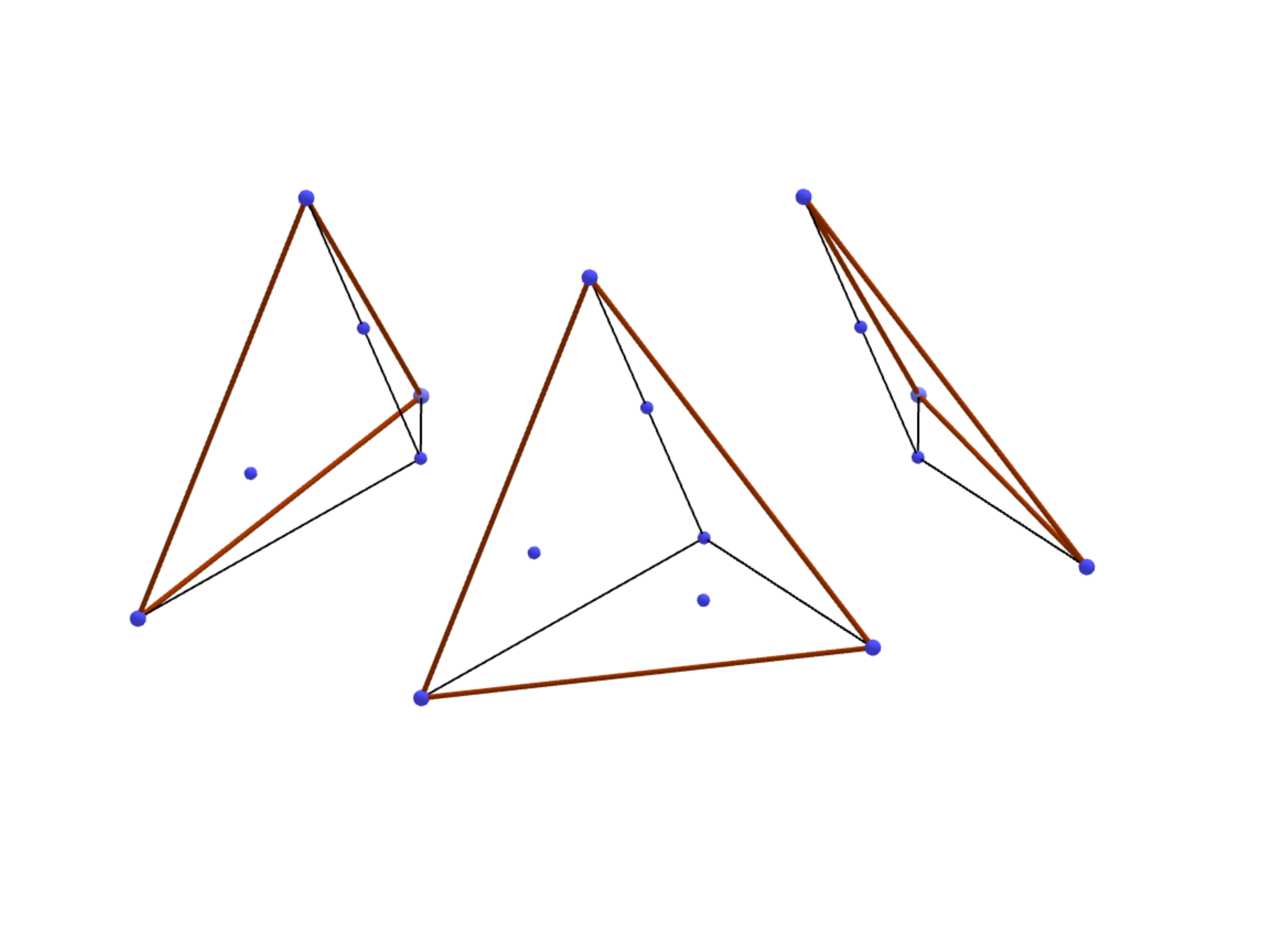}
  \end{center}
  \caption{The first triangulation of $\mathfrak{s}_G$}
  \label{1248-2}
 \end{minipage}
\end{figure}
\end{ex}

The following example shows an application of Proposition \ref{Prop3}.
\begin{ex}\upshape
Let $G=\left\langle\frac{1}{15}(1,6,4,4)\right\rangle$, then the remainder polynomial of $\frac{(1,6,4,4)}{15}$ is 

\begin{eqnarray*}
\mathcal{R}_*\left(\frac{(1,6,4,4)}{15}\right)&=&\frac{(1,6,4,4)}{15}+\frac{(1,3,4,4)}{6}x_2+\frac{(1,2,1,0)}{4}x_3+\frac{(1,2,0,1)}{4}x_4 \\
&+&\frac{(1,0,1,0)}{2}x_3x_2+\frac{(1,0,0,1)}{2}x_4x_2.
\end{eqnarray*}

The coefficient $\frac{1}{6}(1,3,4,4)$ is the coefficient of an iterated term of which the age equals to $2$. By Proposition \ref{Prop3}, $\compl4/G$ has no crepant resolutions.
\end{ex}

Next example shows that the converse proposition to Theorem\ref{thm1} does not hold.
\begin{ex}\upshape
Let $G=\left\langle\frac{1}{24}(1,5,6,12)\right\rangle$. The remainder polynomial of $\frac{(1,5,6,12)}{24}$ is as follows,

\begin{eqnarray*}
\mathcal{R}_*\left(\frac{(1,5,6,12)}{24}\right)&=&\frac{(1,5,6,12)}{24}+\frac{(1,1,1,2)}{5}x_2+\frac{(1,5,0,0)}{6}x_3+\frac{(1,5,6,0)}{12}x_4 \\
&+&\frac{(1,1,1,1)}{2}x_2x_4+\frac{(1,3,1,0)}{5}x_4x_2+\frac{(1,5,0,0)}{6}x_4x_3.
\end{eqnarray*}

The age of the coefficient $\frac{1}{2}(1,1,1,1)$ is $2$, but it is not of an iterated term. Therefore, Theorem \ref{thm1} and  Proposition \ref{Prop3} can not be applied to this case.  This means that $\CFF_{{\bm e}_1}(\sigma)$ does not give a crepant resolution.
However, $\compl4/G$ has a crepant resolution. Actually, Fig.\ref{15612} shows the triangulation of the junior simplex $\mathfrak{s}_G$ corresponding to a crepant resolution.

\begin{figure}[H]
  \begin{center}
   \includegraphics[width=100mm]{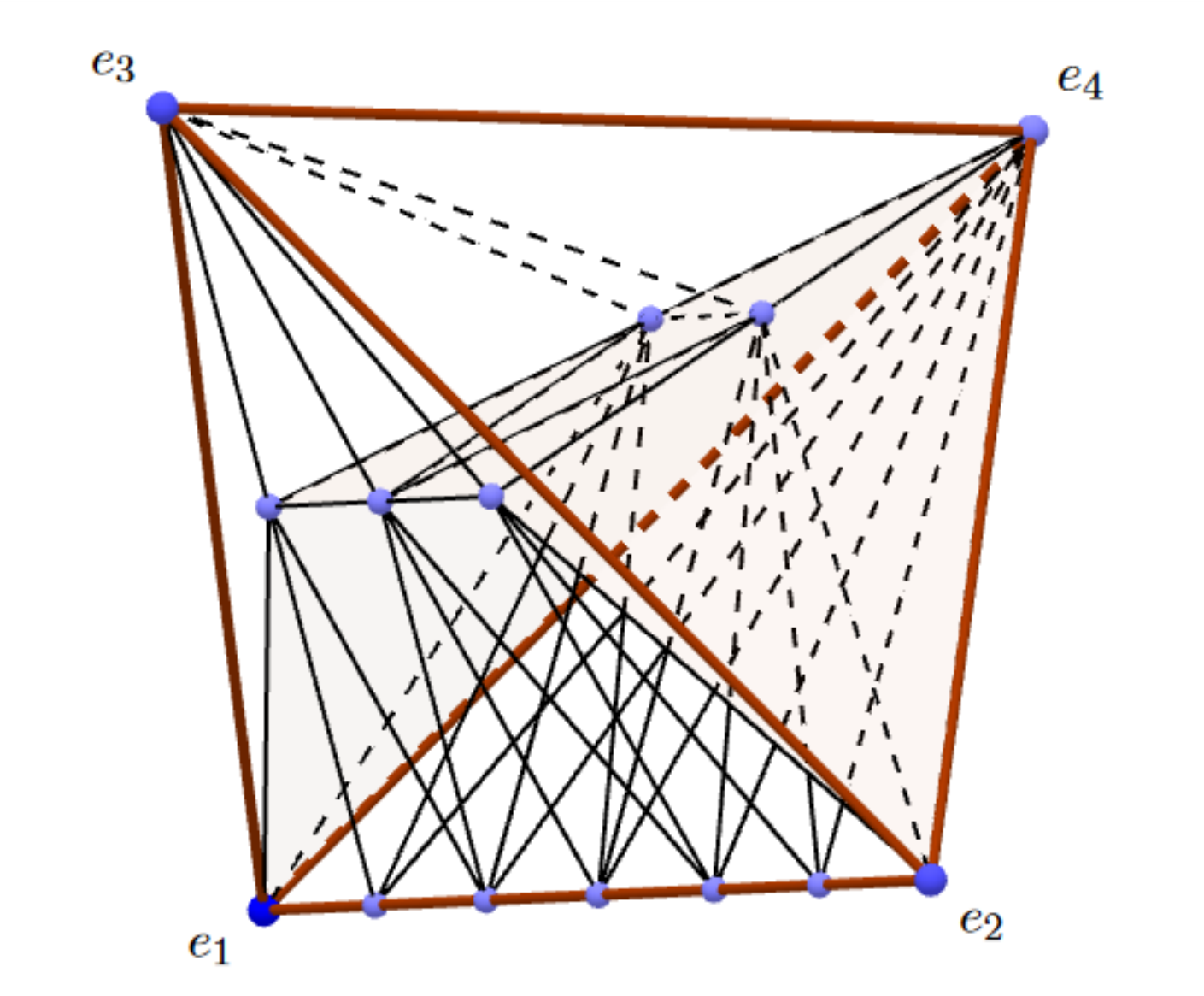}
  \end{center}
  \caption{The triangulation of $\mathfrak{s}_G$ corresponding to a crepant resolution}
  \label{15612}
\end{figure}

\end{ex}

\section{Iterated Fujiki-Oka resolutions for Gorenstein abelian quotient singularities}\label{Ab}

In this subsection, we give a way to construct Fujiki-Oka resolutions for Gorenstein abelian quotient singularities by using Ashikaga's continued fractions repeatedly. As the goal of this chapter, we prove that iterated Fujiki-Oka resolutions for three dimensional Gorenstein abelian quotient singularities are crepant.

\subsection{Basic Generating Systems of $G$}

Let $G\subset\slmc{n}$ be a finite abelian subgroup. Since all the elements in $G$ are simultaneously diagonalizable, there exists a conjugacy class of $G$ which is generated by diagonal matrices. Therefore, we may  assume that $G$ is generated by diagonal matrices. By Proposition \ref{age-prop}, if $G\subset \slmc{3}$, then it is possible to take elements in $G$ of which age is one as the generators of $G$. In higher dimensional case, we assume that the ages of all generators of $G$ are one, because it is clear that $\compl{n}/G$ has no crepant resolutions if the ages of a generator $g$ and the inverse $g^{-1}$ in $G$ are more than one by Proposition \ref{Prop1}. Therefore, we assume the ages of the generators of $G$ are one. By the fundamental theorem of finite abelian groups and the Chinese remainder theorem, there exist a generating system of $G$ as follows:
\[ \left\{ \frac{1}{r_1}(a_{11}, a_{12},\ldots, a_{1n}),\frac{1}{r_2}(0, a_{22},\ldots, a_{2n}),\ldots, \frac{1}{r_{n-1}}(0,\ldots,0, a_{n-1\ n-1},a_{n-1\ n}) \right\} \]
where $r_i, a_{ij}\ (1\leq i\leq n-1,\ i \leq j \leq n)$ are positive integers satisfying $\lcm(r_1,\ldots,r_{n-1})=|G|$ and the following conditions:
\begin{itemize}
\item[(i)] if $a_{ii}=0$, then $a_{ij}=0$ for $i\leq j \leq n$,
\item[(ii)] if $a_{ii}\neq0$, then $a_{ii}=1$ and $\displaystyle \sum_{j=i}^{n}a_{ij} =r_i$.
\end{itemize}
In this paper, we call generating systems of $G$ satisfying the above conditions a {\it basic generating systems} of $G$. Additionally, $G$ can be decomposed to the cyclic components as follows:

\[ G\cong  \left\langle \frac{1}{r_1}(a_{11}, a_{12}, a_{13})\right\rangle \times \cdots \times \left\langle\frac{1}{r_{n-1}}(0,\ldots,0, a_{n-1\ n-1},a_{n-1\ n}) \right\rangle.\]

Clearly, every cyclic component can be decomposed to the product of $p$-Sylow subgroups.

\subsection{Iterated Fujiki-Oka resolutions}

We shall introduce the {\it iterated Fujiki-Oka resolutions} in general dimension. Let $G\subset \slmc{n}$ be a finite abelian subgroup and $H$ be a component of a decomposition by cyclic subgroups of $G$. If the singularity $\compl{n}/H$ is semi-isolated, then we have the Fujiki-Oka resolution $(\widetilde{Y_H}, \mathrm{FO}_1)$ and the toric partial resolution $(Y_G,\phi)$ satisfying the following diagram:

\[
\xymatrix@R=6pt{
                              &&&  \compl{n}\ar[dd]^{\pi_H} \\ 
                            && &                 \\ 
                     \widetilde{Y_H} \ar[rrr]^(.5){\mathrm{FO}_1}_{\textrm{ Fujiki-Oka\ Resolution}}\ar[dd]_{\pi_{G/H}}       &&& \compl{n}/H \ar[dd]^{\pi_{G/H}}  \\ 
                          && \circlearrowright \qquad &  \\ 
                     \widetilde{Y_H}\Big/(G/H) = Y_G \ar[rrr]^{\phi}_{{\rm\qquad  Toric\ Partial\ Resolution}}  &&& \compl{n}/G                  }
\]
where $\pi_{H}$ (resp. $\pi_{G/H}$) is the quotient map by $H$ (resp. $G/H$). Let $X(N_G,\Sigma_{\phi})=Y_G$. If all maximal cones in $\Sigma_{\phi}$ are semi-unimodular with respect to $N_G$, then we have a Fujiki-Oka resolutions $(\widetilde{Y_G},\mathrm{FO}_2)$ for the quotient singularities corresponding to the maximal cones in $\Sigma_{\phi}$.
\[
\xymatrix@R=6pt{
\widetilde{Y_G} \ar[rrrr]^{\mathrm{FO}_2}_{\textrm{ Fujiki-Oka\ Resolution}}   &&&&  Y_H\Big/(G/H) }\]

We note that every singularity in $Y_G$ corresponding to a maximal cone in $\Sigma_{\phi}$ is at worst Gorenstein cyclic quotient singularity which is canonical but not terminal because of the construction.
\begin{defi}\label{it.FO}\upshape
We call the resolution $(\widetilde{Y_G}, \mathrm{FO}_2\circ \phi)$ in the above diagrams an {\it iterated Fujiki-Oka resolution} of $\compl{n}/G$.
\end{defi}

Let $G^{\prime}$ be a finite abelian subgroup which acts on $\widetilde{Y_G}$ equivariant with the torus action and $G$ be a component of a decomposition by cyclic subgroups of $G^{\prime}$. Let $Y_{G^{\prime}}=X(N_{G^{\prime}}, \Sigma_{\phi^{\prime}})$. If all maximal cones in $\Sigma_{\phi}$ are again semi-unimodular with respect to $N_{G^{\prime}}$, then we have a new iterated Fujiki-Oka resolution by extending the above diagram.

\[
\xymatrix@R=10pt{
&                                                                                          &                                                                                               &\compl{n}\ar[dd]^{\pi_H}                       \\ 
&                                                                                          &                                                                                               &                                                       \\ 
&                                                                                          &\widetilde{Y_H} \ar[r]^{\mathrm{FO}_1}\ar[dd]_{\pi_{G/H}}                   &\compl{n}/H \ar[dd]^{\pi_{G/H}}             \\ 
&                                                                                          &                                                                                               &                                                       \\ 
&\widetilde{Y_G} \ar[r]^{\mathrm{FO}_2}\ar[dd]_{\pi_{G^{\prime}/G}}  &Y_G \ar[r]^{\phi}_{\rm T.P.R.}&\compl{n}/G \ar[dd]^{\pi_{G^{\prime}/G}}& \\
                                                   &                                    &                                                                                          &                                                                                               &                                                       \\ 
\widetilde{Y_{G^{\prime}}} \ar[r]^{\mathrm{FO}_3}&Y_{G^{\prime}} \ar[rr]^{\phi^{\prime}}_{{\rm T.P.R.}}                       &                                                                                                &\compl{n}/G^{\prime}                        \\
}
\]
As $(\widetilde{Y_{G^{\prime}}},\mathrm{FO}_3\circ \phi^{\prime})$ in the above, iterated Fujiki-Oka resolutions can be extended under the suitable conditions. We also call these resolutions and the ordinary Fujiki-Oka resolutions {\it iterated Fujiki-Oka resolutions}. 

\begin{lem}\label{lem4-1}
Let $G\subset \slmc{n}$ be a finite abelian subgroup. There exist at least one iterated Fujiki-Oka resolution for $\compl{n}/G$.
\end{lem}
\proofname
Let $\left\{ \frac{1}{r_1}(a_{11}, a_{12},\ldots, a_{1n}),\ldots, \frac{1}{r_{n-1}}(0,\ldots,0, a_{n-1\ n-1},a_{n-1\ n}) \right\}$ be a basic generating system of $G$. we set 
\[H_1=\left\langle \frac{1}{r_{n-1}}(0,\ldots,0, a_{n-1\ n-1},a_{n\ n}) \right\rangle.\]
Then, we have the Fujiki-Oka resolution $X(N_1, \Sigma_1)$ of the singularity $\compl{n}/H_1$ such that the maximal cones in $\Sigma_1$ are obtained from subdividing the two dimensional junior simplex $\mathfrak{s}_2$ spanned by ${\bm e}_{n-1}$ and ${\bm e}_{n}$ into $r_{n-1}$ equal sections. Let $E_i$ be the edge of which endpoints are $\frac{i-1}{r_{n-1}}{\bm e}_{n-1} + \frac{r_{n-1}-i+1}{r_{n-1}}{\bm e}_{n}$ and $\frac{i}{r_{n-1}}{\bm e}_{n-1} + \frac{r_{n-1}-i}{r_{n-1}}{\bm e}_{n}$ for $i=1,\ldots,r_{n-1}$ on $\mathfrak{s}_2$.

As the next step, we set
\[H_2=\left\langle \frac{1}{r_{n-1}}(0,\ldots,0,a_{n-2\ n-2}, a_{n-2\ n-1},a_{n-2\ n}) \right\rangle \times \left\langle \frac{1}{r_{n-1}}(0,\ldots,0, a_{n-1\ n-1},a_{n-1\ n}) \right\rangle.\]
We have the quotient map $\pi_{H_2/H_1}: \compl{n}/H_1 \to \compl{n}/H_2 =X(N_2,\Sigma_1)$. Focus the three dimensional junior simplex $\mathfrak{s}_3$ spanned by ${\bm e}_{n-2}, {\bm e}_{n-1}$ and ${\bm e}_n$. By the definition of basic generating system, there are no lattice points on the edges $E_i\subset \mathfrak{s}_2\subset \mathfrak{s}_3$ for all $i$. Therefore, every maximal cone in $\Sigma_1$ is semi-unimodular, and we have an iterated Fujiki-Oka resolution $X(N_2,\Sigma_2)$.

By repeating similar operation to the above for the subgroup sequence:
\[ H_1\subset H_2 \subset \cdots \subset H_{n-1}=G, \]
we have the sequence of iterated Fujiki-Oka resolutions:
\[ \widetilde{Y_{H_1}}=X(N_1,\Sigma_1), \widetilde{Y_{H_2}}=X(N_2,\Sigma_2), \ldots, \widetilde{Y_G}=X(N_{n-1},\Sigma_{n-1}). \]

\qed\\

By applying Theorem \ref{thm1}, Proposition \ref{Prop2} and Lemma \ref{lem4-1} to the iterated Fujiki-Oka resolutions, we have the following theorem.

\begin{thm}\label{thm2}
Let $\widetilde{Y_{H_1}}, \widetilde{Y_{H_2}}, \ldots, \widetilde{Y_{H_k}}=\widetilde{Y_G}$ be the sequence of iterated Fujiki-Oka resolutions for an $n$-dimensional Gorenstein abelian quotient singularity $\compl{n}/G$. If the ages of all the coefficients in the remainder polynomials associated with every $\widetilde{Y_{H_i}}\ (i=1,\ldots, k)$ are $1$, then the corresponding iterated Fujiki-Oka resolution $\widetilde{Y_G}$ for $\compl{n}/G$ is crepant. 
\end{thm}

Theorem \ref{thm2} and Corollary \ref{maincor} lead to the following corollary.

\begin{cor}\label{cor2}
Assume that $G$ is a finite abelian subgroup of $\slmc{3}$.
Then a crepant iterated Fujiki-Oka resolution exists for $\compl3/G$.
\end{cor}

\subsection{Examples of Iterated Fujiki-Oka resolutions}

At first, we shall see an example of iterated Fujiki-Oka resolutions in three dimension.

\begin{ex}\upshape
Let $G=\left\langle \frac{1}{4}(1,3,0),\frac{1}{4}(1,0,3)\right\rangle$, then $X=\compl3/G$ has a Gorenstein hypersurface singularity defined by $xyz-w^4=0$.
In this case, we have the set $\{ \frac{1}{4}(1,2,1), \frac{1}{4}(0,3,1) \}$ as a basic generating system of $G$. 
According to Lemma \ref{lem4-1}, we set $H=\left\langle\frac{1}{4}(0,3,1)\right\rangle$.
Then the junior simplex of the iterated Fujiki-Oka resolution is transformed as Fig. $6$. 

\begin{figure}[H]\label{ex-itFO1}
  \begin{center}
   \includegraphics[width=120mm]{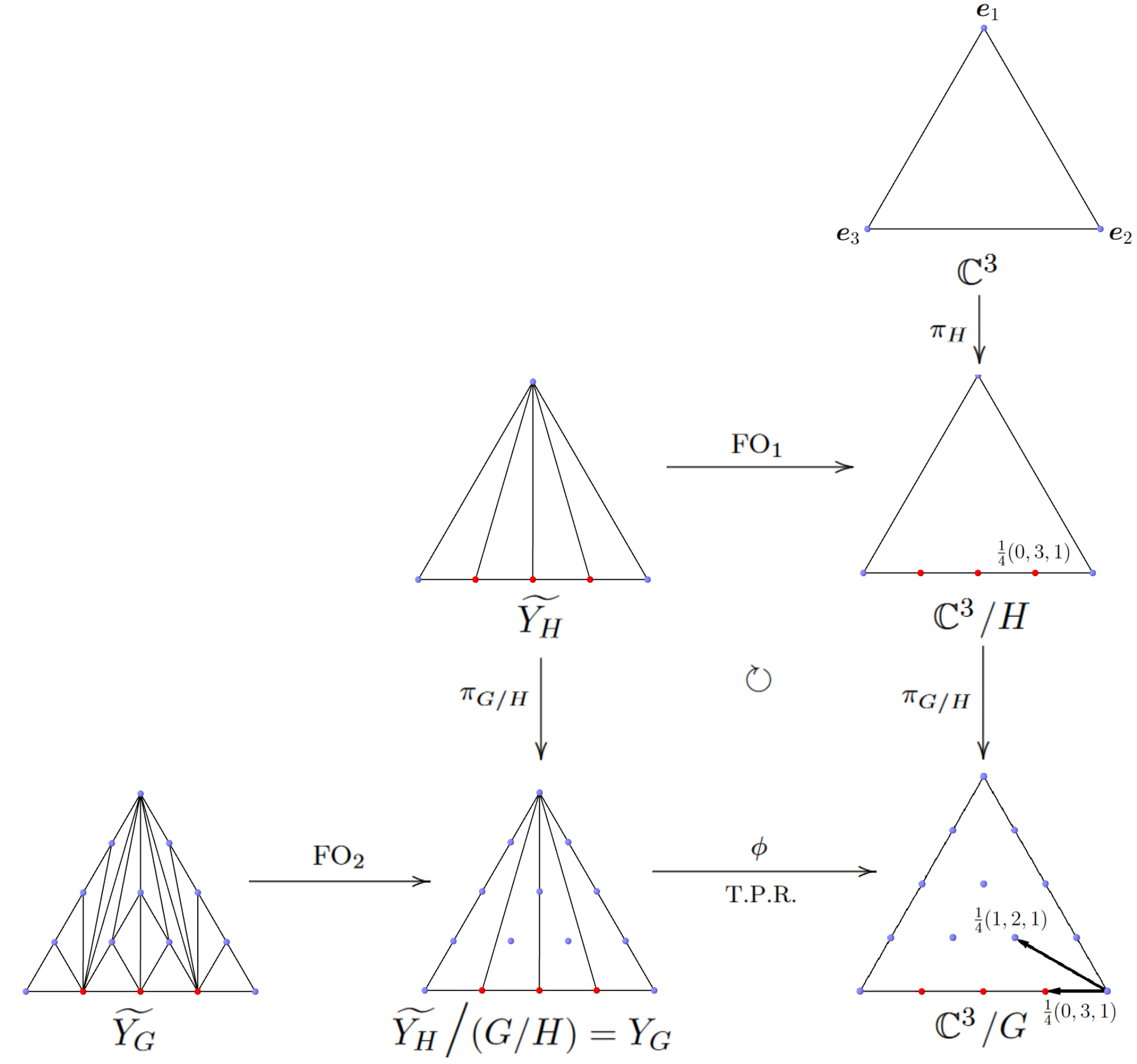}
  \end{center}
  \caption{The iterated Fujiki-Oka resolution for $\left\langle \frac{1}{4}(1,3,0),\frac{1}{4}(1,0,3)\right\rangle$}
\end{figure}
On the other hand, if we choose $\frac{1}{4}(1,2,1)$ as a generator instead of $\frac{1}{4}(0,3,1)$, then we obtain an iterated Fujiki-Oka resolution via a subgroup $H'=\left\langle\frac{1}{4}(1,2,1)\right\rangle$ (see Fig. $7$). In general, the iterated Fujiki-Oka resolution is not unique, and it depends on the choice of the generator.

\begin{figure}[H]\label{ex-itFO2}
  \begin{center}
   \includegraphics[width=120mm]{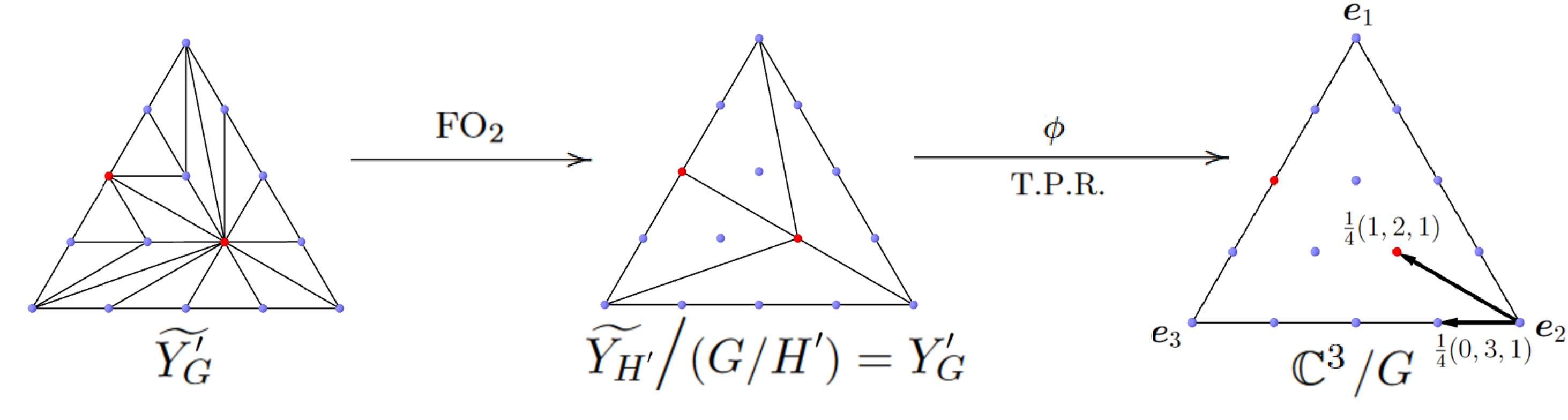}
  \end{center}
  \caption{The iterated Fujiki-Oka resolution via $H'=\left\langle \frac{1}{4}(1,2,1)\right\rangle$}
\end{figure}

\end{ex}

The next example is in four dimensional case.

\begin{ex}\upshape
Let $G=\langle \frac{1}{2}(1,1,0,0), \frac{1}{2}(1,0,1,0), \frac{1}{2}(1,0,0,1)\rangle$, then $X=\compl4/G$ has Gorenstein canonical hypersurface singularity. It is known that $X$ has crepant resolutions. However, $G\mathchar`-\Hilb(\compl4)$ is not a crepant resolution, it is a blow-up of certain crepant resolutions.\\
We can obtain a crepant resolution of $X$ by iterated Fujiki-Oka resolutions. Let $H=\left\langle\frac{1}{2}(1,1,0,0)\right\rangle \subset G$. In addition, this crepant resolution is not blow-down of $G\mathchar`-\Hilb(\compl4)$. In general, the iterated Fujiki-Oka resolution and $G\mathchar`-\Hilb$ give different fans. 

\begin{figure}[H]
  \begin{center}
   \includegraphics[width=100mm]{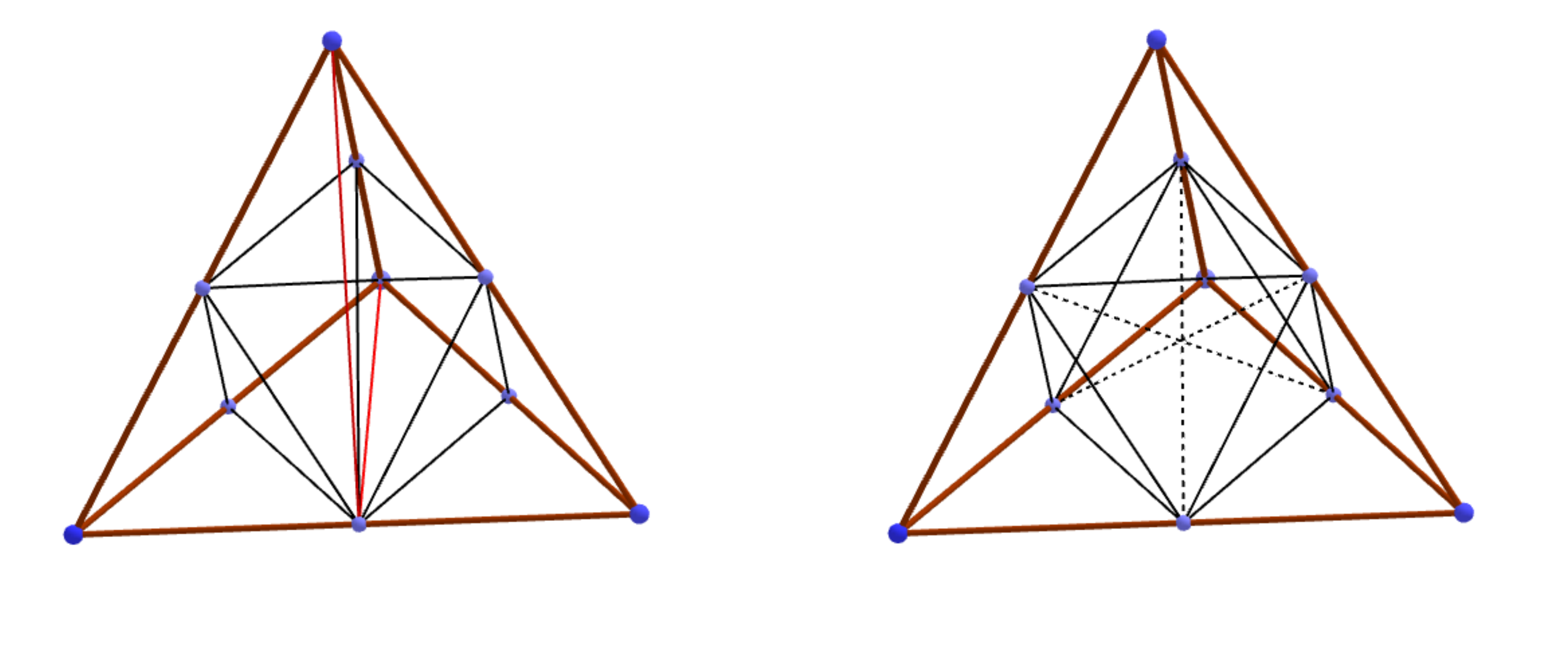}
  \end{center}
  \caption{The junior simplex of the iterated Fujiki-Oka resolution(left) and $G\mathchar`-\Hilb$ (right)}
  \label{HilbandCFF}
\end{figure}
\end{ex}


\section{Discussion}\label{Discuss}
In this section, we shall discuss relationship between the Fujiki-Oka resolution and $\AHilb$ for finite abelian subgroups of $\glmc3$. 
Let's start with a simple case. If $G=\langle\frac{1}{r}(1,1,r-2)\rangle \subset \slmc{3}$, then $\compl{3}/G$ admits a unique projective crepant resolution which coincides with $\AHilb$. On the other hand, the Fujiki-Oka resolution of $\compl{3}/G$ is projective, and Corollary \ref{maincor} claims that it is crepant resolution.
So in this case, the Fujiki-Oka resolution $X(N^{\prime},\CFF_{{\bm e}_1}(\sigma))$ is isomorphic to $\AHilb$. 

We tried to find out a condition such that $X(N^{\prime},\CFF_{{\bm e}_1}(\sigma))$ is isomorphic to $\AHilb$, and propose a conjecture on the relationship between $X(N^{\prime},\CFF_{{\bm e}_1}(\sigma))$ and $\AHilb$.
\begin{conj}
Let $G=\left\langle\frac{1}{r}(a,b,c)\right\rangle \subset \glmc{3}$ and  $\compl3/G$ has semi-isolated singularity. If $X(N^{\prime},\CFF_{{\bm e}_1}(\sigma))$ is isomorphic to $\AHilb$, then the projective toric crepant resolution of $\compl3/G$ exists uniquely up to isomorphism.
\end{conj}

According to "How to Calculate $\AHilb$" \cite{4}, $\AHilb$ can be obtained from Hirzebruch-Jung continued fractions at ${\bm e}_1$, ${\bm e}_2$ and ${\bm e}_3$. On the other hand, $\CFF_{{\bm e}_1}(\sigma)$ can be constructed by Hirzebruch-Jung continued fractions at only two points ${\bm e}_2$ and ${\bm e}_3$. Especially, for an isolated Gorenstein quotient singularity $\compl3/G$, continued fraction fan can be obtained from three ways: $\CFF_{{\bm e}_1}(\sigma)$, $\CFF_{{\bm e}_2}(\sigma)$ and $\CFF_{{\bm e}_3}(\sigma)$. 
Thus, it can be shown that $X(N^{\prime},\CFF(\sigma))$ is isomorphic to $\AHilb$ only if these three continued fractional fans coincides with each other. So we believe the above conjecture is true.\\

At the end of this paper, let us mention a relationship between economic resolutions and Fujiki-Oka resolutions. 
In the case that $G=\left\langle\frac{1}{r}(1,a,r-a) \right\rangle \subset \glmc3$, the quotient singularity $X=\compl3/G$ is terminal and has no crepant resolutions. However, $X$ has an economic resolution (See \cite{29}).

\begin{defi}\upshape
Let $G=\left\langle\frac{1}{r}(1,a,r-a)\right\rangle$ and $N'=\inte2+\frac{1}{r}(1,a,r-a)\inte{}$. Let ${\bm v}_i=\frac{1}{r}(i,\overline{ai}^r,\overline{r-ai}^r) \in N'$ for each integer $1\leq i \leq r-1$. The economic resolution of $X=\compl3/G$ is obtained by the consecutive weighted blow-ups at ${\bm v}_1, {\bm v}_2, \dots, {\bm v}_{r-1}$ from 
$X=\compl3/G$
\end{defi}

\begin{prop}{\rm (\cite[p.381]{29})}
Let $f:Y \to X$ be the economic resolution of $X=\compl3/G$. 

\begin{itemize}
\item[(i)] The toric variety $Y$ is smooth and projective over $X$,
\item[(ii)] The morphism $f$ satisfies
$$
K_Y=f^*(K_X)+\sum_{1\leq i \leq r-1}\frac{i}{\;r\;}\;E_i
$$
\end{itemize}

where $E_i$ is the prime exceptional divisor.
\end{prop}

Since the weighted blow-up with ${\bm v}_1, \dots , {\bm v}_{r-1}$ coincides with the Fujiki-Oka resolution $X(N^{\prime},\CFF(\sigma))$, $X(N^{\prime},\CFF(\sigma))$ is an economic resolution. As O. Kedzierski\cite{22} and S. J. Jung\cite{17} showed Theorem \ref{ER}, economic resolutions can be expressed in some moduli spaces. Thus, Fujiki-Oka resolutions can be written as a moduli space in the case of $G=\left\langle\frac{1}{r}(1,a,r-a)\right\rangle$.

\begin{thm}\label{ER}
The economic resolution $Y$ of a three fold terminal quotient singularity $X=\compl3/G$ is isomorphic to the birational component $Y_{\theta}$ of the moduli space $\mathcal{M}_{\theta}$ of $\theta$-stable $G$-constellations for a suitable parameter $\theta$.
\end{thm}

In particular, Kedzierski\cite{19} has shown that $\AHilb$ is an economic resolution in some special cases.

\begin{thm}{\rm (\cite[Section 2.]{19})}
Let $G \subset \glmc3$ be the finite subgroups generated by $\frac{1}{r}(1,a,r-a)$ with $a=1$ or $r-1$. Then $\AHilb$ is isomorphic to the economic resolution of the quotient variety $\compl3/G$.
\end{thm}


Kohei Sato\\
\textsc{National Institute Of Technology, Oyama Collage, 771 Nakakuki, Oyama, Tochigi, 323-0806, Japan.}\\
E-mail:k-sato@oyama-ct.ac.jp\\

\noindent Yusuke Sato\\
\textsc{Graduate School Of Mathematical Sciences, University Of Tokyo 3-8-1 Komaba, Meguro-ku, Tokyo 153-8914, Japan.}\\
E-mail:yusuke.sato@ipmu.jp

\end{document}